\documentclass[twoside,11pt]{amsart}
\usepackage{graphicx, lscape}
\usepackage{amsmath, amsfonts, amssymb, ifthen}
\usepackage{caption,xcolor}
\usepackage{subcaption}
\usepackage{mathrsfs}
\usepackage[justification=centering]{caption}
\usepackage[colorlinks=true]{hyperref}

\newcommand{\CC}{\mathbb{C}}
\newcommand{\DD}{\mathbb{D}}

\newcommand{\EE}{\mathcal{E}}
\newcommand{\cE}{\mathcal{E}}
\newcommand{\cM}{\mathcal{M}}
\newcommand{\cH}{\mathcal{H}}
\newcommand{\cN}{\mathcal{N}}
\newcommand{\cS}{\mathcal{S}}
\newcommand{\sA}{\mathscr{A}}
\newcommand{\Chat}{\hat{\mathbb{C}}}

\newcommand{\Dplus}{D_+}
\newcommand{\Dminus}{D_-}
\newcommand{\Rnba}{R_{n,a,b}}
\newcommand{\Rna}{r_{n,a}}
\newcommand{\mapRna}{r}
\newcommand{\Mn}{{M}_n}
\newcommand{\crpt}{w}
\newcommand{\spine}{\mathcal{S}_n}
\newcommand{\Dmke}{D^{m,k}_{2\ell}}
\newcommand{\Dmko}{D^{m,k}_{2\ell-1}}

\newtheorem{theorem}{Theorem}[section]
\newtheorem{corollary}[theorem]{Corollary}

\newtheorem{conjecture}{Conjecture}

\theoremstyle{definition}
\newtheorem{definition}[theorem]{Definition}
\newtheorem*{notation}{Notation}

\newcommand{\ep}{\varepsilon}

\DeclareMathOperator\Arg{Arg}

\newtheorem{prop}[theorem]{Proposition}
\newtheorem{cor}[theorem]{Corollary}
\newtheorem{lem}[theorem]{Lemma}
\newtheorem{defn}[theorem]{Definition}

\begin{document}

\title[Baby Julia Sets in Parameter Space for rational maps]
{Exploring baby Julia sets in parameter space slices for Generalized McMullen Maps}

\author[S.~Boyd]{Suzanne Boyd}
\address{Department of Mathematical Sciences\\
University of Wisconsin Milwaukee\\
PO Box 413\\
Milwaukee, WI 53201, USA\\
ORCID: 0000-0002-9480-4848}
\email{sboyd@uwm.edu}

\author[K.~Brouwer]{Kelsey Brouwer}
\address{Mathematics Department\\
Aquinas College\\
1700 Fulton St E\\
Grand Rapids, MI 49506\\
USA}
\email{brouwerk@aquinas.edu}

\author[M.~Hoeppener]{Matthew Hoeppner}
\address{Mathematics and Statistics Department\\
Hope College\\
27 Graves Place\\
Holland, MI 49423, USA}
\email{hoeppner@hope.edu}

\date{\today}

\begin{abstract} 
 For the family of complex rational functions of the form 
\noindent \mbox{$\Rnba(z) = z^n + \dfrac{a}{z^n}+b$,} 
 known as ``Generalized McMullen maps'', for $a\neq 0$ and $n \geq 3$ fixed, we 
 describe the apparent phenomena of baby Julia sets in parameter space appearing both in slices with independent critical orbits and a slice defined by imposing a critical orbit relation.  
 Specifically, we introduce the subfamily where one of two critical orbits is set to be a super-attracting fixed point, provide some general results on this subfamily and describe how Julia set copies in the parameter space slice occur--due to parameters for which the other critical orbit is in the (not immediate) basin of attraction of this fixed critical point.  We provide several conjectures on this intriguing phenomena to catalyze further study.
 \end{abstract}

\maketitle

\markboth{\textsc{S.\ Boyd, K.\ Brouwer and M.\ Hoeppner}}
  {\textit{Baby Julia Sets in parameter space for rational maps}}


\section{Introduction: the big picture} 
\label{sec:introduction}

In 1985, Adrien Douady and John Hamal Hubbard (\cite{DouadyHubbard}) revealed and explained the appearance of the, even then already famous, {\em Mandelbrot set} in varied families of iterative processes. The \text{Mandelbrot Set} $\cM$ is the set of $c$-values in the complex plane such that the orbit of $0$ under the map $P_c(z)=z^2+c$ is bounded, so it is a ``boundedness locus''. Its boundary is the family $P_c$'s bifurcation locus. The origin is significant as the sole critical point in this family, and it is of multiplicity two.  The generalization $z\mapsto z^n+c$ is the family with the sole critical point of higher multiplicity.

Douady and Hubbard showed that multiple homeomorphic copies of the Mandelbrot set occur in the bifurcation locus generated from applying Newton's Method to a cubic polynomial family with a single parameter. They defined what it means for a map to behave like $P_c$, calling such a map {\em polynomial-like of degree two}.  We call homeomorphic copies of the Mandelbrot set $\cM$ ``baby'' Mandelbrot sets or ``baby $\cM$'s''.

\begin{figure}
\centering
\begin{subfigure}{0.4\textwidth}
\includegraphics[width=.95\textwidth,keepaspectratio]{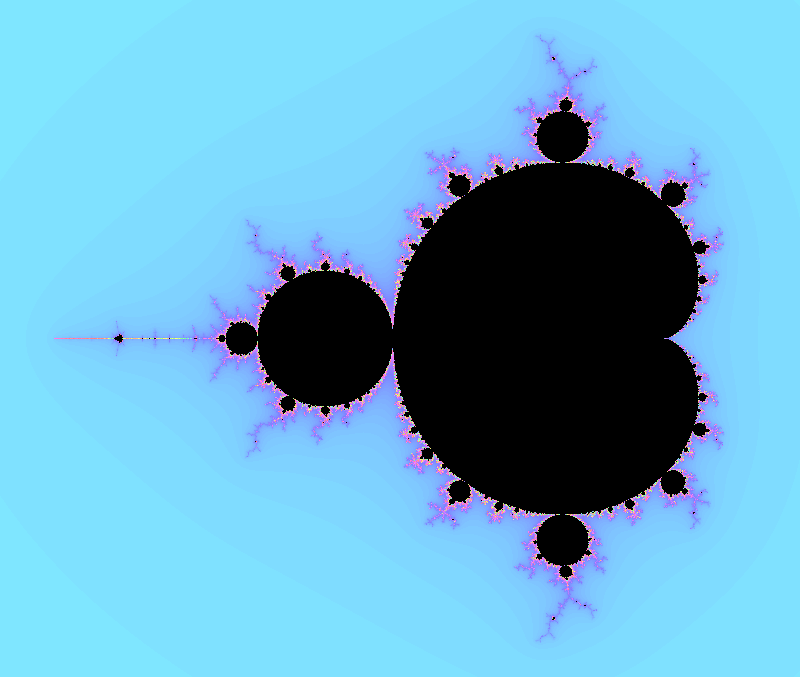}
  \caption{The Mandelbrot Set: the bifurcation locus for $P_c(z)=z^2+c$.}
\end{subfigure}%
\begin{subfigure}{0.5\textwidth}
\includegraphics[width=1.0\textwidth,keepaspectratio]{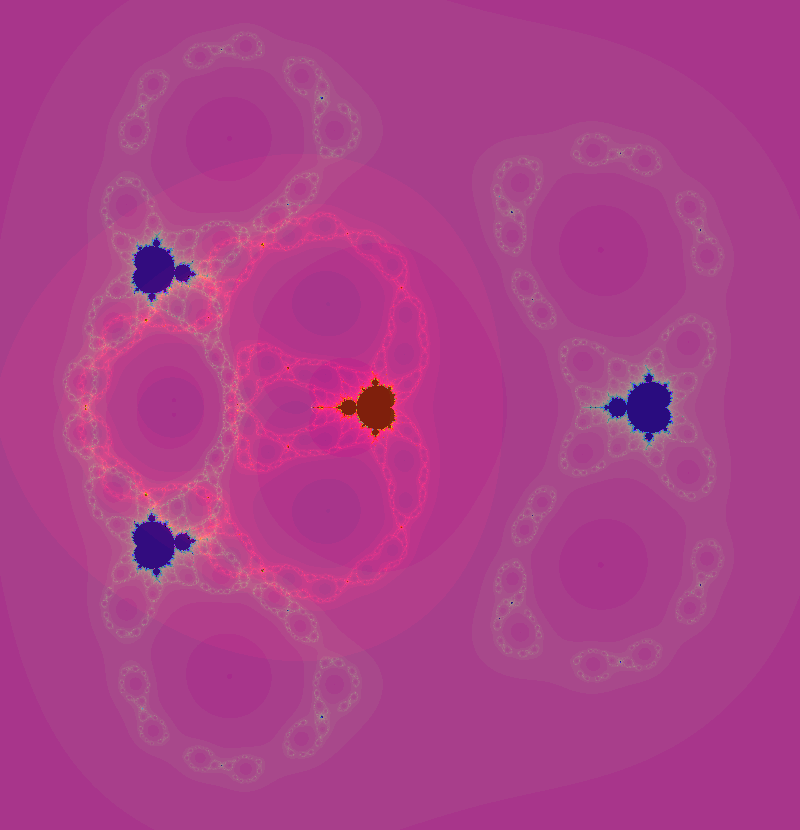}
\caption{The bifurcation locus in the $a$-plane for $R_{n,a,b}(z)=z^n+a/z^n +b,$ with $n=5, b=0.5$.}
\end{subfigure}
\caption{\label{fig:Mandelbrotone}
Parameter planes: for $P_c$ (left), and baby Mandelbrot sets in a slice of a bifurcation locus (right).}
\end{figure}

The Mandelbrot set is on the left of Figure \ref{fig:Mandelbrotone}; the right illustrates apparent small copies of the Mandelbrot set in the bifurcation locus of our more general family of rational maps of interest.

Fifteen years later, Curt McMullen (\cite{mcmullen}) showed the presence of the Mandelbrot set in other families of maps was universal: every non-empty bifurcation locus of any analytic family will contain (quasi-conformal) copies of the Mandelbrot Set, or of the boundedness locus of $z\mapsto z^n+c$ based on the degree of the critical points. 

Douady and Hubbard's original approach prescribes criteria for proving when Mandelbrot set copies exist in specific locations of families of map (which we describe in more detail later in this article). 

As is typical, our approach to exploring a parameter space for a family of maps begins with studying the dynamical plane of maps in this family, including the behavior of critical orbits of the map and the invariant set of interest: the Julia set.  
We first define the \textit{Fatou set} of $P_c$ as the set of $z$-values in the domain with stable behavior: where the iterates of $P_c$ is a normal family in the sense of Montel. The \textit{Julia set} $J$ is the complement to the Fatou set. The \textit{filled Julia set} $K$ is the union of the Julia set and the bounded Fatou components. Then $K$ is the set of all points with bounded orbits, and $J=\partial K$.

A rational family which generalizes $z\mapsto z^n+c$, and in which it is easy to find baby Mandelbrot sets, is the ``Generalized McMullen Family'':
$$
\Rnba(z) = z^n + \dfrac{a}{z^n}+b~,
$$
where $n \geq 3$ is an integer, $a\neq 0$ is complex, and $b$ is complex. Thus, this family's parameter space is effectively $\CC^2.$ 
Note the point at $\infty$ is super-attracting, and $0$ is a preimage of $\infty$. Thus we can still define the filled Julia set $K$ as the set of points with bounded orbits, and get $J=\partial K$.
McMullen introduced the study of this family in the case $b=0$ (\cite{mcmullen_example}), Devaney and colleagues have studied this family in that case and for other generalizations (\cite{DevaneyRussell,DevaneyHalos,DevaneyCase2,DevaneySurvey2013,
DevLimBeh,DevCheckerboard,DevBlanchard,DevBlanchard2,
devgar,devlook,DevRussellConnectivity,
DevaneyKozma}), and the authors of the present article and colleagues have studied the generalized McMullen maps (\cite{BoydBrouwer1,BoydHoeppner1,BoydMitchell,BoydSchulz}) as have Xiao, Qui, and Yongchen (\cite{xiaoqiu}).

Now instead of just one critical point, Generalized McMullen Maps have $2n$ finite critical points: the solutions of $w^{2n}=a$, the $2n^{th}$ roots of $a$. While this seems like a lot of critical orbits to keep track of,  every one of these critical points maps $2:1$ onto one of just two critical values, $v_{\pm} = b \pm 2\sqrt{a}$. Hence, the behavior inherently has two degrees of freedom. 
Figure~\ref{fig:Mandelbrotone} (right) shows a slice of the bifurcation locus, the $a$-plane for $n=5, b=0.5$. To draw this image, for each parameter we iterate $v_-$ and $v_+$, and assign the average of two color values: both critical values yield black if the orbit appears bounded, but one critical value yields blue upon escape, and the other red. So the color assigned is essentially either purple (both escape), dark red (one escapes), dark blue (the other escapes), or black (both appear bounded). 

\begin{figure}
\centering
\begin{subfigure}{.5\textwidth}
  \centering
  \includegraphics[width=.95\textwidth]{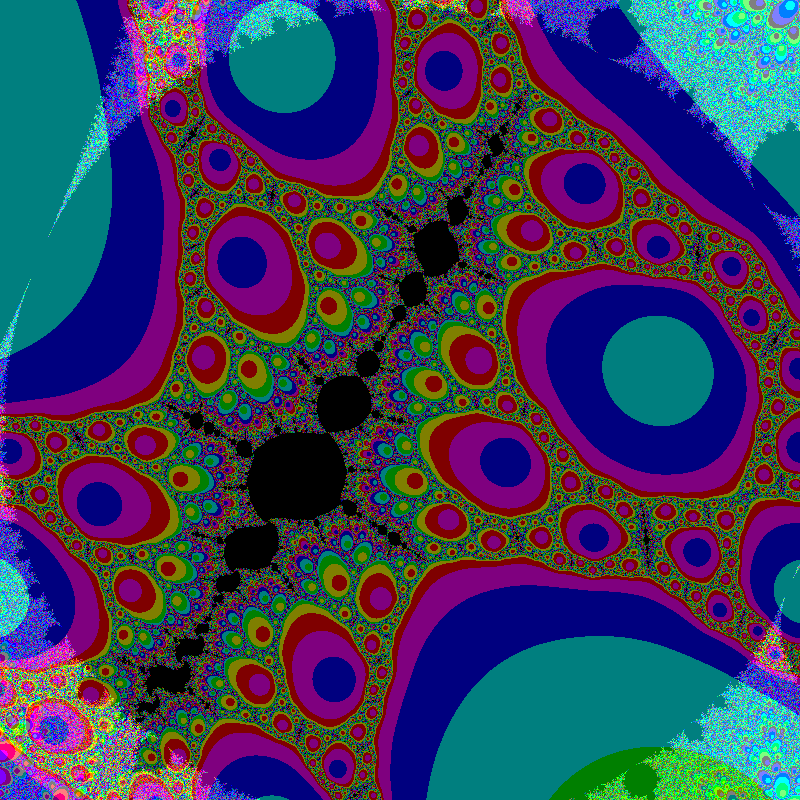}
  \caption{Figure in the necklace of $R_{4,6a,a}$}
  \label{fig:NecklaceZoom2}
\end{subfigure}%
\begin{subfigure}{.5\textwidth}
  \centering
  \includegraphics[width=.95\textwidth]{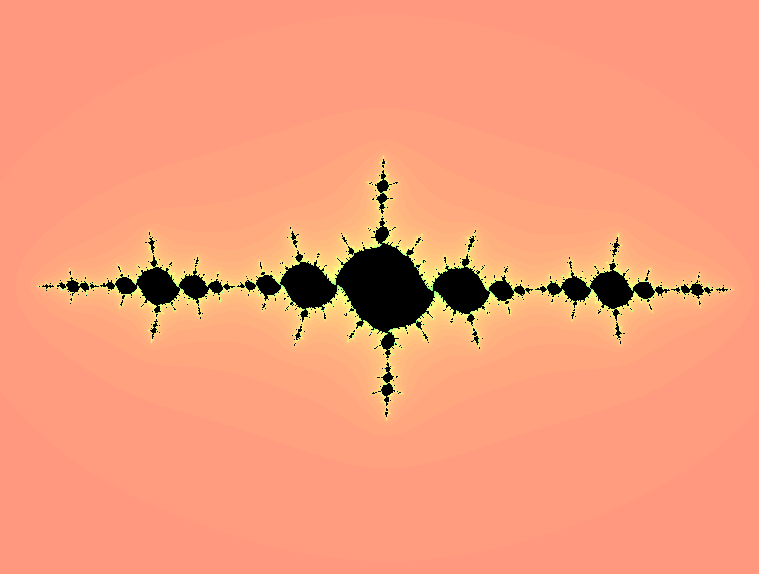}
  \caption{Julia set for $P_c$ with $c=-1.3 + 0.02i$}
  \label{fig:BasicJulia2}
\end{subfigure}%
\caption{Apparent (partial) baby Julia set in a parameter slice of $\Rnba$, here the $a$-plane for $b:=6a$, and the corresponding quadratic Julia set (basillicas in the basillica)}
\label{fig:JuliaComparison2}
\end{figure}

In this article, we reveal and begin the study of a striking behavior, that turns out to be the result of having two independent, bounded critical orbits.
In generating parameter space images for this family, it becomes almost immediately clear there are Mandelbrot-like shapes present. What is not  as readily apparent is that in many cases there are also smaller shapes, buried in the ``necklace structure'' (the chains of shaded loops visible in Figure~\ref{fig:Mandelbrotone}, right), which upon magnification resemble distorted quadratic {\em Julia} sets. In slices where the two critical orbits are independent, thus each has a bifurcation locus (such as, one parameter is held constant while the other is free, or the family with $b:=ta$ for a chosen $t$, introduced in \cite{BoydHoeppner1}), they seem to appear where the bifurcation loci of the two critical points intersect in a specific way.
Indeed,
it is understood (and no surprise given McMullen's universality result) that in 1-dimensional parameter slices of the bifurcation locus of the Generalized McMullen family, in regions where one critical orbit escapes, zooming into the necklace structure of the bifurcation locus for the other critical orbit reveals baby Mandelbrot sets (presumably, densely). On the other hand, we find baby Julia sets in parameter space in a similar way: 

\begin{conjecture}[Location of baby Julia sets in parameter slices of Generalized McMullen Maps]
\label{conj:locationbabyJs}
    In 1-dimensional subfamilies of the bifurcation locus of the generalized McMullen family in which the two critical orbits have independent behavior, where one critical value generates a baby Mandelbrot, if the necklace structure of the bifurcation locus of the other critical point intersects the interior of that baby Mandelbrot, then, rather than baby Mandelbrot sets dense in the necklace structure, there are ``baby'' topological Julia sets, or partial baby Julia sets (partial as in Figure~\ref{fig:JuliaComparison2}). We do not expect these baby Julia sets in parameter slices to be quasi-conformal copies, only homeomorphic.  
\end{conjecture}

See Figure~\ref{fig:JuliaComparison2}. In fact, we have no reason to believe this phenomena is restricted to Generalized McMullen Maps, they are simply the rational family we have explored thoroughly. For this family, 
we have observed this to occur in many computer-generated images of bifurcation locus slices; for example in the $b$-plane for small $a$. In this article we collect some machinery which we hope can be used in a future study establishing this conjecture.

Buff and Henriksen (\cite{BuffHenriksen}) prove examples of Julia sets in cubic parameter space slices defined by a fixed indifferent critical point, though in this article we cannot use their techniques, our bifurcation locus is not simply connected.

Our conjecture above is apparent to the explorer generating the right images. We next provide a more analytical conjecture about conditions under which, in general, topological baby Julia sets appear in the bifurcation locus of the Generalized McMullen subfamilies.
     Our conjecture is an alteration of the Douady and Hubbard criteria for having a baby Mandelbrot in parameter space. Essentially, for baby Mandelbrot sets to occur in parameter space, there must be a sole critical point in a region ($U'_a$) in dynamical space in which the map is degree two polynomial like (from $U'_a$ to its image $U_a$), and a region ($W$) in parameter space such that as the parameter ($a$) loops around the boundary of that region, the critical value loops around the critical point (or, the critical value loops around $U'_a\setminus U_a$) in each dynamical plane (and the regions need to be defined nicely and vary nicely, see Section~\ref{sec:prelim} for more details).   On the other hand:

\begin{conjecture}[Criteria for baby Julia sets in parameter space]
\label{conj:BabyJcriteria}
   Suppose there is a region $W$ in a one-dimensional subfamily of the generalized McMullen parameter space, parameterized by the complex variable $a$, in which one critical value, say $v_+(a)$, is topologically stable; i.e., $W$ is entirely contained inside of a hyperbolic component of the bifurcation locus of $v_+(a)$. Moreover, suppose in $W$ the orbit of $v_+(a)$ is bounded and in fact satisfies the criteria so that $v_+(a)$ lies within a baby polynomial Julia set, within its larger Julia set: there is a region $U'_a$ containing that baby Julia set, and $f_a:U'_a\to U_a$ is ``polynomial-like''.
    
    Now, suppose there is another critical value, say $v_-(a)$, which for some parameter inside of $W$ maps (eventually) onto $v_+(a)$, i.e., $v_-(a)$ lies in a preimage copy of the baby Julia set containing $v_+(a)$.  
    If as $a$ loops around $\partial W$, $v_-(a)$, loops around the corresponding (eventual) preimage of $U'_a \setminus U_a$,  then inside of $W$, in parameter space, the bifurcation locus of $v_-$ sweeps out a homeomorphic copy of the baby polynomial Julia set associated with $v_+(a)$.    
\end{conjecture}

 We expect proof of this second conjecture is needed to establish the first. 

Now that we have set the general stage, we focus in on a 1-dimensional subfamily of the generalized McMullen family, for which we provide concrete results.  Even without the independent critical orbit behavior described in Conjecture~\ref{conj:locationbabyJs}, we can find apparent baby Julia sets in the bifurcation locus in parameter space slices defined by certain critical orbit relations. For example, for the bulk of the remainder of this paper, we study the subfamily defined by requiring the critical value $v_+$ to be a super-attracting fixed point, which takes the form of a one-dimensional complex manifold in the $(a,b)$-space. 
Specifically, we set $b$ so that $v_+$ is equal to the principle critical point which is the root of $a^{1/2n}$ of argument $\Arg(a)/2n$, and consider the subfamily:
\begin{equation}
\label{eqn:subfamilyrna}
\Rna(z) = z^n + \frac{a}{z^n}+ (a^{1/2n}-2\sqrt{a}), 
\end{equation}
where the $a^{1/2n}$ is the principle root,
the degree $2n$ is always fixed, and $n\geq 3$. 
Now, the other critical value $v_-$ is free. 
Since $v_+$ is a fixed point, the behavior is determined entirely by the critical orbit of $v_-$. 

\begin{figure}
\centering
\begin{subfigure}{0.5\textwidth}
  \centering
\includegraphics[width=.95\textwidth,keepaspectratio]{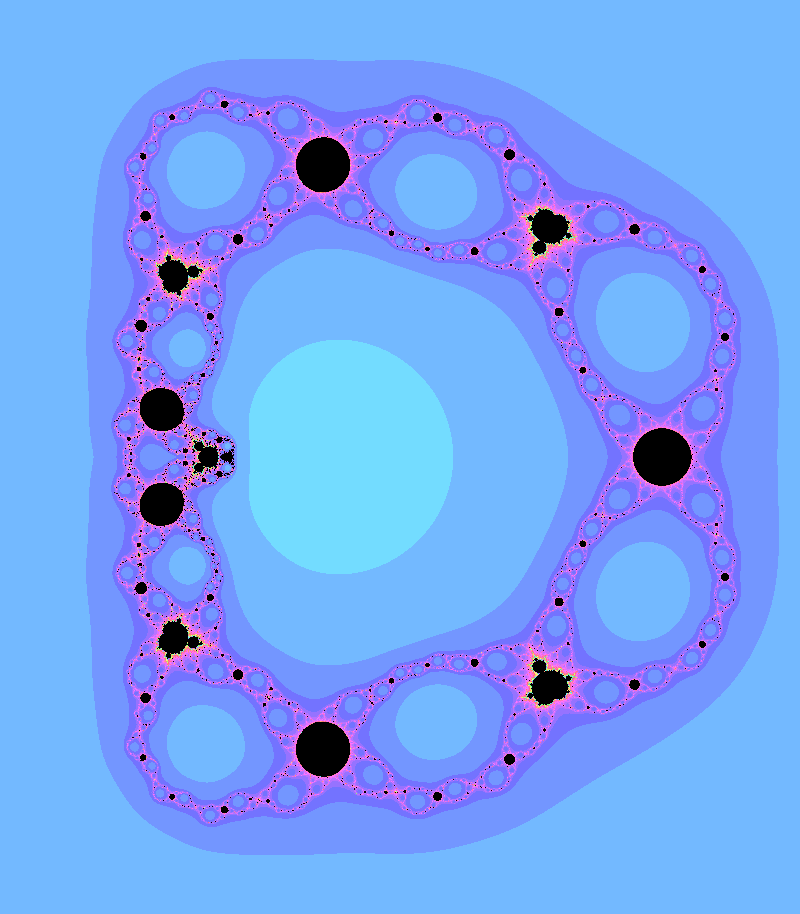}
  \caption{Parameter plane for $n=6$}
  \label{fixed_crit_n=6_aplane_intro}
\end{subfigure}%
\begin{subfigure}{0.5\textwidth}
  \centering
\includegraphics[width=.95\textwidth,keepaspectratio]{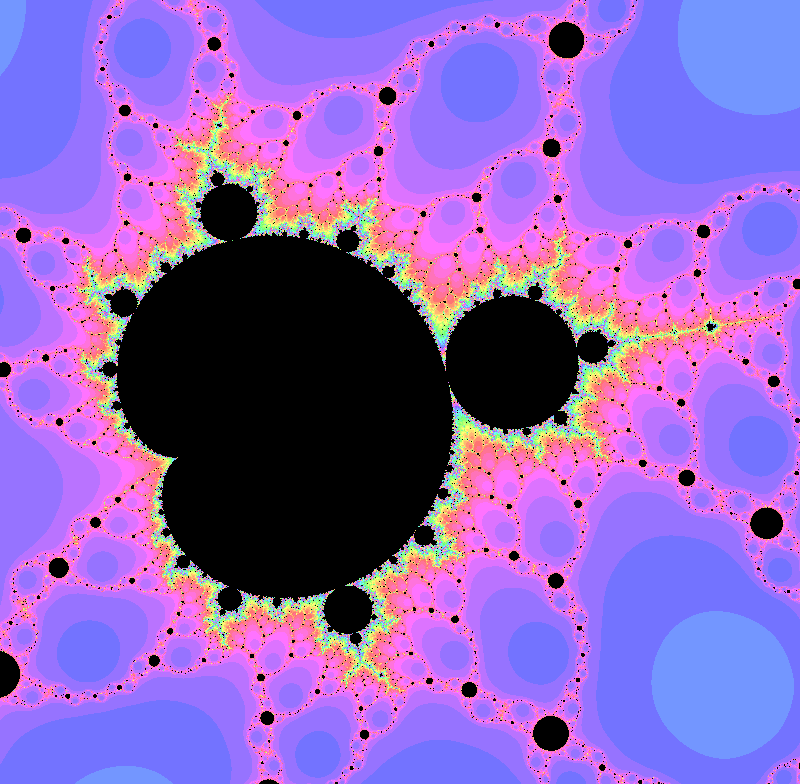}
  \caption{Zoomed in on baby Mandelbrot set}
  \label{fixed_crit_n=6_zoom_intro}
\end{subfigure}%
\caption{The parameter slice of $R_{6,a,b}$ with one critical value held at a fixed point.}
\label{fig:Fixed_Crit_Point_intro}
\end{figure}

We observe in slices of this type that not only are there regions strewn throughout the slice in which the behavior of $v_-$ is essentially unconstrained and generates a baby $\cM$, the full range of quadratic polynomial behavior, but also throughout the parameter plane there appear to be topological disks. See Figure~\ref{fig:Fixed_Crit_Point_intro} for an example of the bifurcation locus for $n=6$.  

\begin{figure}
\centering
\begin{subfigure}{.5\textwidth}
  \centering
  \includegraphics[width=.95\textwidth,keepaspectratio]{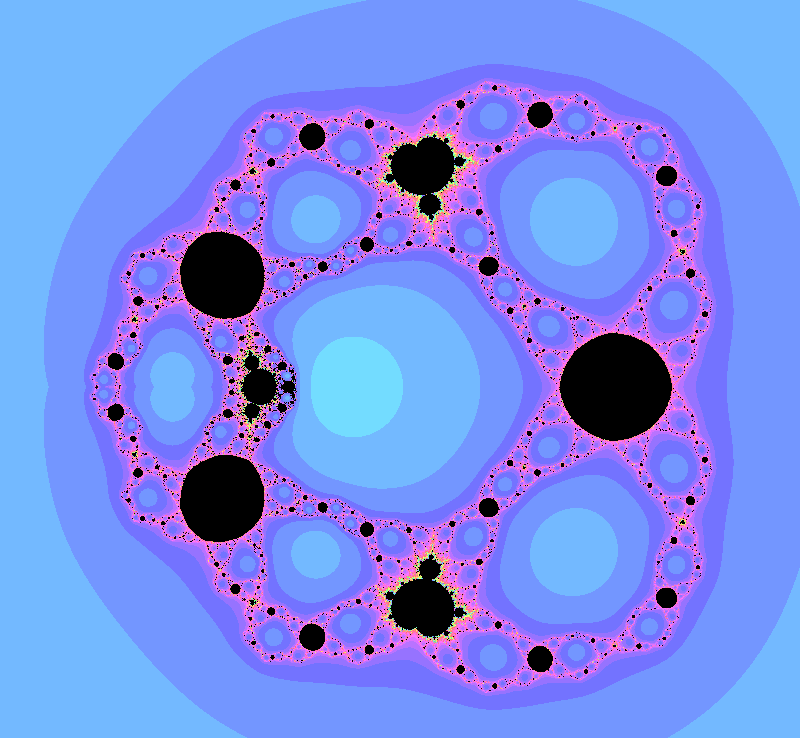}
  \caption{Parameter plane for $n=4$}
  \label{fixed_crit_n=4_aplane}
\end{subfigure}%
\begin{subfigure}{.5\textwidth}
  \centering
  \includegraphics[width=.95\textwidth,keepaspectratio]{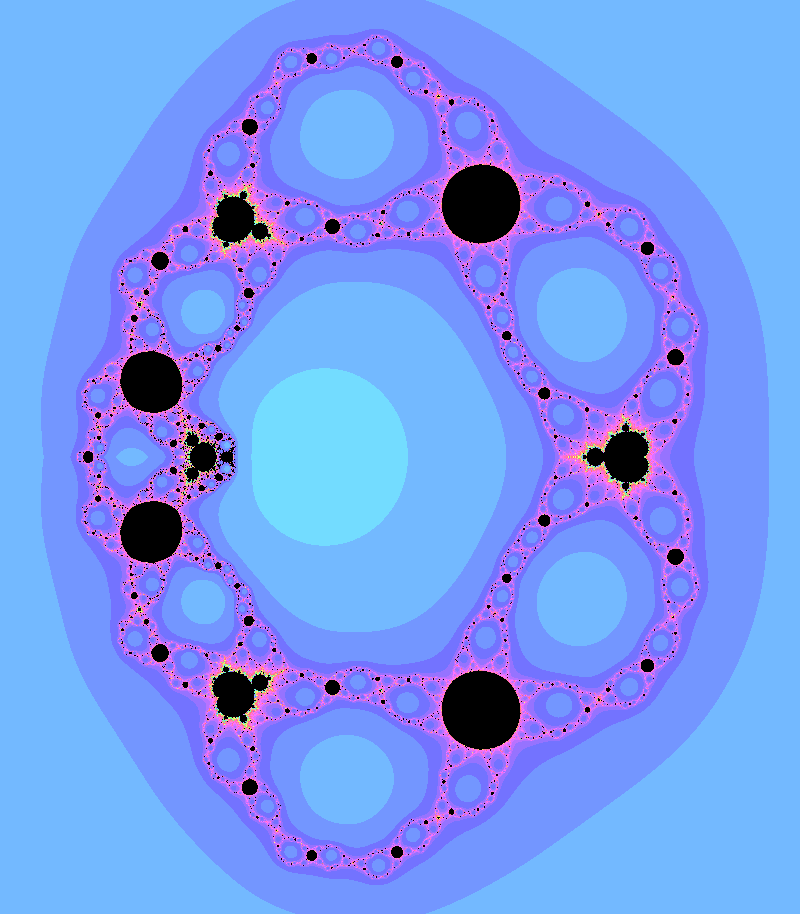}
  \caption{Parameter plane for $n=5$}
  \label{fixed_crit_n=5_aplane}
\end{subfigure}
\begin{subfigure}{.5\textwidth}
  \centering
    \includegraphics[width=.95\textwidth,keepaspectratio]{Mfixedcrit_n_6.png}
  \caption{Parameter plane for $n=6$}
  \label{fixed_crit_n=6_aplane2}
\end{subfigure}%
\begin{subfigure}{.5\textwidth}
  \centering
  \includegraphics[width=.95\textwidth,keepaspectratio]{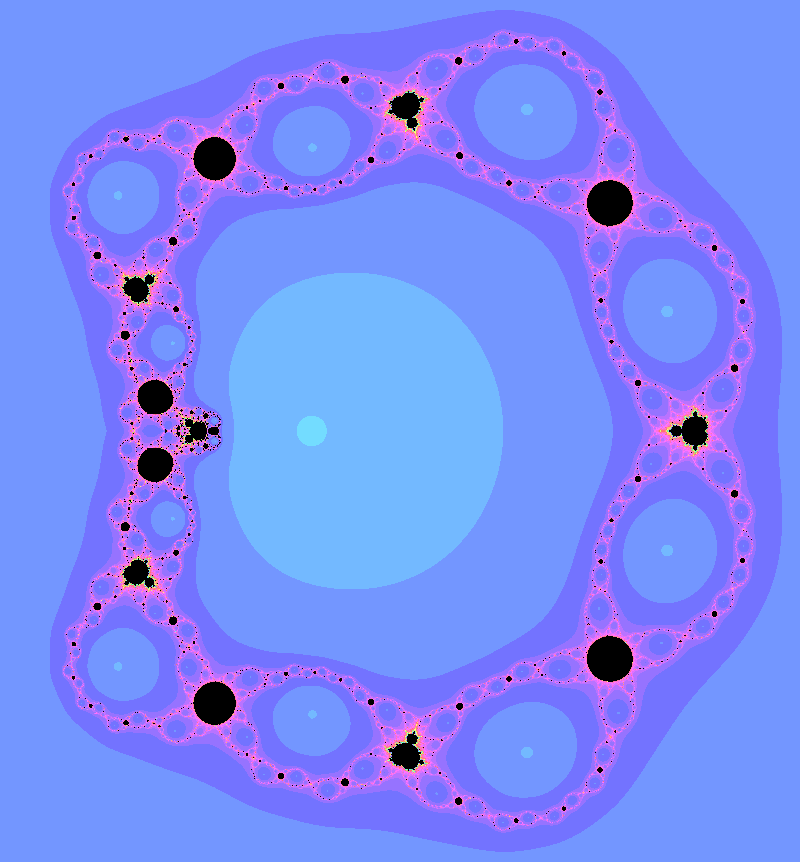}
  \caption{Parameter plane for $n=7$}
  \label{fixed_crit_n=7_aplane}
\end{subfigure}
\caption{Parameter planes for $\Rna$ for $n=4,5,6,7$.}
\label{fig:FixedCritPointExamples1}
\end{figure}

As $n$ increases we see some interesting patterns emerge, see Figure \ref{fig:FixedCritPointExamples1}. There are what look like $n-2$ well-defined baby $\cM$'s in each plane, as well as $n-1$ large filled topological disks. There also appear to be more disks and $\cM$'s throughout the necklace structure in these images, but these $n-1$ disks and $n-2$ baby $\cM$'s are the most prominent, so we'll refer to them as the ``principal'' copies. Near $a=0$, there is also a shape present in each which almost looks like two baby $\cM$'s colliding, thus suggesting the dynamics in that area is not the same as in an isolated baby $\cM$. 

Now, these are not mere topological disks; $P_0(z)=z^2$ is the quadratic polynomial map with one fixed critical point (at $0$) and thus seems related to our situation, and its Julia set is a disk (indeed, Julia sets in this map's hyperbolic component are all topological disks). 

Our main result for this subfamily is
a description of these topological disks in the parameter space, in terms of the dynamics the polynomial $p_0(z)=z^2$, whose Julia set is the unit disk $\DD$.

\begin{theorem}
    \label{thm:hypcompdisk}
The bifurcation locus in the subfamily $r_{n,a}(z)= z^n + \frac{a}{z^n}+ (a^{1/2n}-2\sqrt{a}),$
contains  
$n-1$
hyperbolic components, $\cH_{2j}$ for $j=1, \ldots, n-1$,
for which there are homeomorphisms
$\Phi_j \colon \cH_{2j} \to \mathbb{D}$. 

Each $\cH_{2j}$ has dynamical center at 
$$
 a_{2j} = \left( \frac{1-e^{i\frac{2j\pi}{n}}}{4} \right)^{\frac{2n}{n-1}},
$$
where at the center, $\mapRna(v_+)=v_+$ is a super-attracting fixed point, and $\mapRna(v_-)=v_+$, and such that $v_-$ is not in the same  Fatou component as $v_+$.
\end{theorem}

To further describe the subfamily $\Rna$, we also provide some results on the location of the boundedness locus, including  Lemma~\ref{lem:spinecritfixed}, in which we calculate a spine for the bifurcation locus in this slice, which appears to approach the cardioid $a=\frac{1}{16}(1+e^{i\theta})^2$) as $n\to \infty$, Lemma~\ref{lem:M_in_annulus}, calculating an annulus centered at $1/8$ containing the boundedness locus for sufficiently large $n$, yielding Corollary~\ref{cor:M_in_disk} that the boundedness locus lies in the disk centered at the origin of radius just over $1/2$ for sufficiently large $n$. Then, in Proposition~\ref{prop:c-patterns}, we calculate $2n-1$ parameters $a_k$ such that for $k$ odd, the map fixes $v_-$, and for $k$ even, $v_-$ maps to $v_+$.  The odd index parameters are conjectural centers of baby Mandelbrot sets. 
For the hyperbolic components containing the even index parameters, we utilize a result from \cite{BoydBrouwer1} to describe the dynamics of each map in one of these components (Proposition~\ref{prop: preimage classes}), and we describe the appearance of the topological disks in parameter space. 

One avenue for further study would be to further describe the qualities of the map $\Phi$, and study its extension to the boundary. 
We expect the construction in the proof of Theorem~\ref{thm:hypcompdisk} is the core idea for the map linking parameter spaces in general with baby Julia sets. Though we do not state any further details in this article, we note another interesting generalization would be to consider how Julia sets appear in the two complex dimensional parameter space, rather than in a one-dimensional slice.

We close this opening section with the organization of the sections. In Section \ref{sec:prelim} we provide background information.
Section~\ref{sec:dpr} contains dynamical plane results pertinant to this subfamily.
In Section~\ref{sec:spines}, we provide some results on the location of the boundedness locus for this subfamily.  
In Section~\ref{sec:centers} we calculate the dynamical centers of the proposed principle baby $\cM$'s and topological disks, proving Proposition~\ref{prop:c-patterns}. In Section~\ref{sec:babyMandels}, we provide a conjecture ($n-2$ baby $\cM$s in the parameter plane for this subfamily) and suggestions on approaching its proof by using results from this article as well as \cite{BoydHoeppner1,BoydMitchell}). 
In section~\ref{sec:paramdisks}, we study the dynamics in the topological disks in parameter space, establish Theorem~\ref{thm:hypcompdisk}. 

\subsection*{Acknowledgements} We thank John Hamal Hubbard and Sarah Koch for helpful conversations and suggestions. We would also like to thank Brian Boyd for the computer program \textit{Dynamics Explorer} used to generate all of the Mandelbrot and Julia images in this paper.

\section{Preliminaries}
\label{sec:prelim}

In this section we describe useful tools and results from earlier work. We also begin with an overview of previous work on this family. 

Generalized McMullen maps, including the one-parameter subfamily with $b=0$, have been studied previously by Devaney and colleagues as well as the first author and colleagues.
In \cite{DevaneyRussell,DevaneyKozma}, Devaney and coauthors study the family in the case of $b$ at the center of a hyperbolic component of the Mandelbrot set for $P_b$ (that is, the critical point is a fixed point). 
 For $n\geq 2$, Devaney and colleagues study the subfamily with $b=0$, ``McMullen maps", in papers such as \cite{DevaneyHalos,DevaneySurvey2013}. 
 They establish the location of $n-1$ baby $\cM$'s in the slice of $\Mn(\Rnba)$ in the $a$-parameter plane when $b=0$ (see Figure \ref{fig:Mandelbrotone} (right) for an example). In the $b=0$ case there is only one free critical orbit, thus the bifurcation locus is the boundary of the set of parameters for which the critical orbit is bounded (the ``boundedness locus'' in this case). 
 In \cite{devgar}, Devaney and Garijo  study Julia sets as the parameter $a$ tends to $0$, for a different generalization of McMullen maps: $z \mapsto z^n + \dfrac{a}{z^d}$ with $n \neq d$.
In \cite{DevaneyHalos,jangso, Stoertz-error} Devaney, Marotta, Stoertz, and co-authors find $n$ baby $\cM$'s when $n\neq d$.

In \cite{BoydSchulz}, the first author and Schulz study the geometric limit as $n\to \infty$ of Julia sets and of $\Mn(\Rnba)$, for $\Rnba$ for any complex $c$ and any complex $a\neq 0$.  
In \cite{BoydMitchell}, the first author and Mitchell establish the location of baby $\cM$'s in $a$-plane slices in the parameter space for $b\neq 0$ and $b\in [-1,1]$ fixed (and $n\geq 3$), and in  $b$-plane slices for $n\geq 1$ and $\frac{1}{10} \leq a \leq 4$.
In \cite{BoydHoeppner1}, the first and last author of the present article establish the location of $n$ baby  $\cM$'s in the $a$-parameter place for fixed $|b|\geq 6$ and $n$ large enough that $4|b|+8 \leq 2^{n+1}$. Additionally, we begin the study the ``linear'' slices with $b=ta,$ for any $t$ fixed, and locate a spine and a neighborhood of the spine in which the bifurcation locus must lie. 

Not all of the behavior of generalized McMullen maps is polynomial-like. For example, McMullen showed there are Julia sets which are Cantor sets of simple closed curves (\cite{devlook,mcmullen_example}). 
Xiao, Qiu, and Yin (\cite{xiaoqiu}) establish a topological description of the Julia sets (and Fatou components) of $\Rnba$ according to the dynamical behavior of the orbits of its free critical points. This work includes a result that if there is a critical component of the filled Julia set which is periodic while the other critical orbit escapes, then the Julia set consists of infinitely many homeomorphic copies of a quadratic Julia set, and uncountably many points. In \cite{BoydBrouwer1}, the first two authors of this article provide a combinatorial model of the dynamics in the case that one critical orbit is in a baby quadratic Julia set with a period two cycle, and the other critical orbit strictly eventually lands inside of that baby. As we shall see, that phenomena is what drives the creation of copies of quadratic Julia sets in parameter space.

First, some notation.

\textbf{Notation.} Let $\cN(S)$ denote a small neighborhood (of unspecified size) of a set $S$. 

\subsection{Polynomial-Like Maps}
Douady and Hubbard \cite{DouadyHubbard} provide criteria for when baby polynomial Julia sets exist in the Julia sets of other maps, and for when baby $\cM$'s exist within the parameter planes of other families, by showing that a particular family of functions behaves locally like a degree two polynomial.

\begin{defn}[\cite{DouadyHubbard}, Chapter 1] \label{defn:DH1}
A \textit{polynomial-like map of degree $d$} is a proper, analytic map of degree $d$, $f:U' \to U$, where $U'$ and $U$ meet the following conditions: both $U'$ and $U$ are open subsets of $\CC$ which are homeomorphic to disks, $U' \subset U$, and $U'$ is relatively compact in $U$.

The \textit{filled Julia set} of a polynomial-like map is the collection of points in $U'$ whose orbits never leave $U'$, which we will denote by $K(f|_{U'})$.
\end{defn}

So note $K(f|_{U'}) = \cap_{n} f^{-n}(U')$, which is likely not the same as $K(f) \cap U'$.

Each polynomial-like map of degree $d$ has exactly $d-1$ critical points. Here we are  only concerned with polynomial-like maps of degree two. A polynomial-like map of degree two is a two-to-one map from $U'$ to $U$, with a unique critical point in $U'$.  The ``straightening theorem'' shows this is an apt name.

\begin{theorem}[\cite{DouadyHubbard}, Chapter 1, Theorem 1] \label{thm:DH2}
A polynomial-like map $f$ of degree $d$ is topologically conjugate on its filled Julia set to some degree $d$ polynomial $P$ on the polynomial's filled Julia set. 

Moreover, for $d=2$, the conjugacy $\phi$ is a ``hybrid-equivalence'': quasi-conformal in  a neighborhood  $\cN( K(f|_{U'}))$ to a neighborhood $\cN(K_P)$, and conformal on $K(f|_{U'})$, so
$$\phi:\cN(K_P) \to \cN(K(f|_{U'}))$$ 
and
$$
P \circ \phi = \phi\circ f.
$$

Finally (again assuming $d=2$), if $K(f|_{U'})$ is connected, then $P$ is unique up to conjugation by an affine map. 
\end{theorem}

We call the filled Julia set $K(f|_{U'})$ of the polynomial-like map a ``baby'' Julia set.

For Generalized McMullen Maps, the degree of the polynomial map is always $d=2$.

\begin{definition}[\cite{DouadyHubbard}, Chapter 2]
\label{defn:DH-familyf}
Let $\Lambda$ be a complex analytic manifold and ${\mathrm{\bf f}} = (f_\lambda: U'_{\lambda} \to U_\lambda )_{\lambda \in \Lambda}$ a family of polynomial-like mappings. Set $\mathscr{U} = \{ (\lambda,z) | z \in U_\lambda \} $, $\mathscr{U'} = \{ (\lambda,z) | z \in U'_\lambda \} $, and $f(\lambda,z) = (\lambda, f_\lambda(z))$. Then ${\mathrm{\bf f}}$ is an \textit{analytic family} if the following hold:
\begin{enumerate}
    \item $\mathscr{U}$ and $\mathscr{U}'$ are homeomorphic over $\Lambda$ to $\Lambda \times \DD$.
    \item The projection from the closure of $\mathscr{U}'$ in $\mathscr{U}$ to $\Lambda$ is proper.
    \item The mapping $f:\mathscr{U}' \to \mathscr{U}$ is complex analytic and proper. 
\end{enumerate}
The degree of each $f_\lambda$ is the same if $\Lambda$ is connected, call this the degree of $\textbf{f}$. 
\end{definition}

\begin{definition}[\cite{DouadyHubbard}, Chapter 2]
\label{defn:DH-Mf}
    Let ${\mathrm{\bf f}}$ be an analytic family of degree $2$ polynomial-like mappings parameterized by $\Lambda$.

    Set $\mathscr{K}_{\textbf{f}}  = \{ (\lambda, z) | z \in K(f_\lambda|_{U'_\lambda}). \}$ 
    Note $\mathscr{K}_{\textbf{f}} $ is closed in $\mathscr{U}$ and the projection of $\mathscr{K}_{\textbf{f}} $ onto $\Lambda$ is proper, since $\mathscr{K}_{\textbf{f}} = \cap_{n} f^{-n}(\mathscr{U'})$.

    Let $\mathrm{M}_{\textbf{f}}$ denote the set of $\lambda$ in $\Lambda$ for which $K(f_\lambda|_{U'_\lambda})$ is connected and
 $\mathscr{R}=\Lambda - \partial \mathrm{M}_{\textbf{f}}$.
\end{definition}

Now, since $d=2$ in our case of interest, by the straightening theorem, we can associate each member $f_\lambda$ of an analytic family of polynomial-like maps of degree $2$ with a polynomial $P_{\chi(\lambda)} = z^2 + \chi(\lambda)$
and with a homeomorphism between neighborhoods of their Julia sets 
$$\phi_\lambda:\cN(K(f_\lambda|_{U'_\lambda})) \to \cN (K(P_{\chi(\lambda)})), $$
where $\cN(K(f_\lambda|_{U'_\lambda}))  \subset U'_\lambda$, and where $\phi_{\lambda}$ is a hybrid equivalence from $f_\lambda$  to $P_{\chi(\lambda)}$ (conformal on filled Julia sets, quasi-conformal in neighborhoods), with
$$
P_{\chi(\lambda)}\circ \phi_\lambda = \phi_\lambda\circ f_\lambda.
$$
Further, Douady and Hubbard investigate how $\phi_\lambda$ and $P_{\chi(\lambda)}$ vary with $\lambda$. The following are their results, again in the case of degree $d=2$.

\begin{theorem}[\cite{DouadyHubbard}, Chapter 2, Theorems 1 and 2]
\label{thm:DHC2T1}
Suppose ${\mathrm{\bf f}} =(f_\lambda : U'_\lambda \to U_\lambda)_{\lambda \in \Lambda}$ 
is an analytic family of polynomial-like mappings of degree $d=2$.

Then each $f_\lambda$ is conjugate via $\phi_\lambda$ to $z \mapsto z^2 + \chi(\lambda),$ 
where the mapping $\chi : \Lambda \to \CC$ is continuous on $\Lambda$ and analytic on $\mathrm{\mathring{{M}}}_{\mathrm{\bf f}}$. 
And further, for  $\lambda \in \mathscr{R}(=\Lambda - \partial \mathrm{M}_{\mathrm{\bf f}})$,
$\phi_\lambda$   depends continuously on $\lambda$. 
\end{theorem}

Douady and Hubbard translate into the polynomial-like map setting the Ma\~ne, Sad, and Sullivan (\cite{MSS}) results  
on structural stability in a neighborhood of the Julia set $\partial K(f|_{U'})$ (the Julia sets, {\em not} the filled Julia sets). 

\begin{prop}[\cite{DouadyHubbard}, Chapter 2, Proposition 10]
\label{prop:DH_tau}
For any $\lambda_0$ not in $\partial \mathrm{M}_{\mathrm{\bf f}}$, there's a map $\tau_\lambda$ defined for $\lambda$ in a neighborhood of $\lambda_0$ (avoiding $\partial \mathrm{M}_{\mathrm{\bf f}}$), which is holomorphic in $\lambda$ and quasi-conformal in $z$ (with a dilitation ratio bounded by a constant independent of $\lambda$), and which conjugates $f_{\lambda_0}$ to $f_\lambda$ in a neighborhood of their respective Julia sets: $\tau_\lambda \circ f_{\lambda_0} = f_\lambda \circ \tau_\lambda$, with $\tau_{\lambda_0} = \text{id}$.
\end{prop}

We can improve upon this using McMullen and Sullivan's \cite{McMSull1998}. In their Chp.\ 7, they describe the structural stability extended to the Riemann sphere, $\Chat$.  Recall a \textit{holomorphic family} of rational maps $f_\lambda$ over a complex manifold $X$ is a holomorphic map $X\times \Chat \to \Chat$, given by $(\lambda, z) \mapsto f_\lambda(z).$
In a holomorphic family of rational maps, $f_{\lambda_0}$ is \textit{topologically stable} if for all $\lambda$ in some neighborhood of $\lambda_0$, $f_\lambda$ and $f_{\lambda_0}$ are topologically conjugate (on the sphere). 
They show this implies there is in fact a quasi-conformal conjugacy, in their Theorem 7.1, and they further show that the set of critical orbit relations is constant. 

The following is the definition of holomorphic motion, as given in \cite{McMSull1998}.
\begin{definition} 
    \label{defn:holmot}
A \textit{holomorphic motion} of a set $S \subset \Chat$, over a connected complex manifold with basepoint $(X, \lambda_0)$, is a mapping $\tau \colon X \times S \to \Chat$, given by 
$(\lambda,z) \mapsto \tau_\lambda(z)$, such that:
\begin{itemize}
    \item For each fixed $z\in S$, $\tau_\lambda (z)$ is a holomorphic function of $\lambda$;
    \item For each fixed $\lambda \in X$, $\tau_\lambda(z)$ is an injective function of $z$; and 
    \item The injection is the identity at the basepoint (that is, $\tau_{\lambda_0}(z)=z.)$
\end{itemize}
\end{definition}  

McMullen-Sullivan's method establishes the following.

\begin{theorem}[\cite{McMSull1998}, Theorem 7.4]
\label{thm:MS_tau}
    For a topologically stable $\lambda_0$, there is a neighborhood of $\lambda_0$ and a {holomorphic motion} $\tau_\lambda$ of the sphere, over $\lambda$ and its neighborhood, which ``respects the dynamics''; that is, 
$\tau_\lambda(f_{\lambda_0}(z)) = f_\lambda(\tau_\lambda(z))$, whenever $z$ and $f_{\lambda_0}(z)$ are both in the domain of $\tau_{\lambda}$.
\end{theorem}

They use the Harmonic $\lambda$-lemma of Bers-Royden (see \cite{SullThur1986} and \cite{BersRoy1986}), which shows a holomorphic motion of a set $A$ in a neighborhood in parameter space can be extended \textbf{uniquely} to a holomorphic motion of the entire sphere, over a smaller neighborhood in parameter space (and the motion and its extension agree on the domain overlap). The uniqueness is what allows them to guarantee the compatibility with the dynamics. 

In our setting, compatibility with the dynamics implies that the $\tau$'s of Proposition~\ref{prop:DH_tau} and Theorem~\ref{thm:MS_tau} coincide, though Douady and Hubbard's $\tau_\lambda$ was only defined near the Julia sets. This is why we used the same symbol for both of these.

Moreover, McMullen-Sullivan apply their results to hyperbolic quadratic polynomials. Let $H\subset \CC$ be the open set of hyperbolic parameters. Partition $H$ by setting $H(\infty)$ to the set of parameters whose critical orbits escape to $\infty$; $H_0(p)$ is the set of parameters whose critical point is periodic of order $p$; and $H(p)$ is the set of parameters whose critical orbit is strictly attracted to a cycle of period $p$, so that $H_0(p)$ is not in $H(p)$. 

Theorem 10.1 (\cite{McMSull1998}) states first that each connected component of $H(p)$ is isomorphic to a punctured disk \textit{and} represents a single quasiconformal conjugacy class. The multiplier of the attracting cycle gives a natural isomorphism of $\DD^*$. Secondly, for each integer $p\geq 1, H_0(p)$ is finite, and in fact $H_0(p)$ are the punctures in the centers of the components in $H(p)$.
In general, one must be careful about centers of hyperbolic components, as $\phi_{\lambda}$ is valid near $J$, and extends via holomorphic motion to all of $\CC$ for parameters \textit{except} hyperbolic component centers.
Thirdly, $H(\infty)$ is also isomorphic to a punctured disk \textit{and} represents a single quasiconformal conjugacy class. Finally, maps in different components of these sets represent different conjugacy classes. 

We use the above results in our setting, especially in Sections~\ref{sec:dpr} and~\ref{sec:paramdisks}.

\medskip

Finally, we give Douady and Hubbard's criteria for the existence of baby Mandelbrot sets in parameter space, as stated by Devaney based on Douady and Hubbard:

\begin{theorem}\cite{DevaneyHalos,DouadyHubbard} \label{thm:DH3}
Suppose we have a family of polynomial-like maps $f_\lambda : U'_\lambda \to U_\lambda$ which satisfy the following:
\begin{enumerate}
    \item The parameter $\lambda$ is contained in an open set in $\CC$ which contains a closed disk $W$;
\item 
 The boundaries of $U'_\lambda$ and $U_\lambda$ both vary analytically as $\lambda$ varies;
\item The map $(\lambda, z) \to f_\lambda(z)$ depends analytically on $\lambda$ and $z$;
\item Each $f_\lambda$ is polynomial-like of degree two and has a unique critical point, $c_\lambda$.
\end{enumerate}
Suppose that for each $\lambda$ in the boundary of $W$ we have that $f_\lambda(c_\lambda) \in U_\lambda \setminus U'_\lambda$, and that $f_\lambda(c_\lambda)$ winds once around $U'_\lambda$ (and therefore once around $c_\lambda$) as $\lambda$ winds once around the boundary of $W$. Then the set of all $\lambda$ for which the orbit of $c_\lambda$ does not escape from $U'_\lambda$ is homeomorphic to the Mandelbrot set.

\end{theorem}

\begin{figure}
  \includegraphics[width=0.45\linewidth,keepaspectratio]{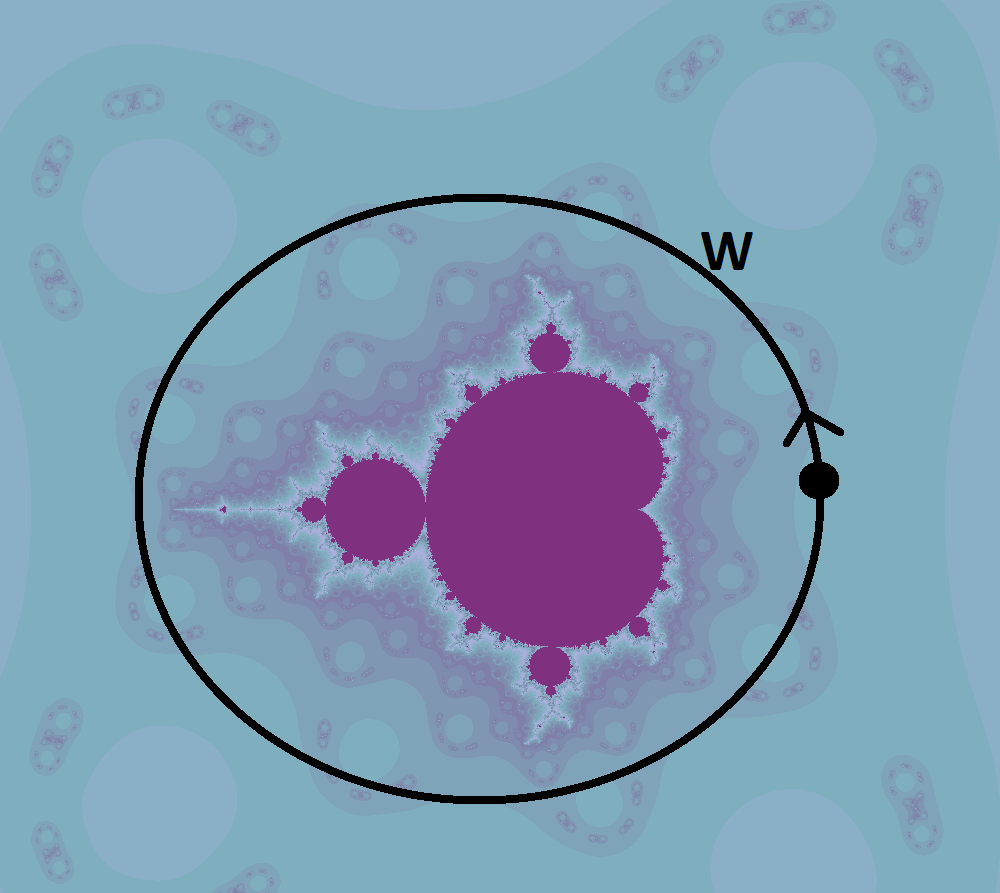} 
  \includegraphics[width=0.45\linewidth,keepaspectratio]{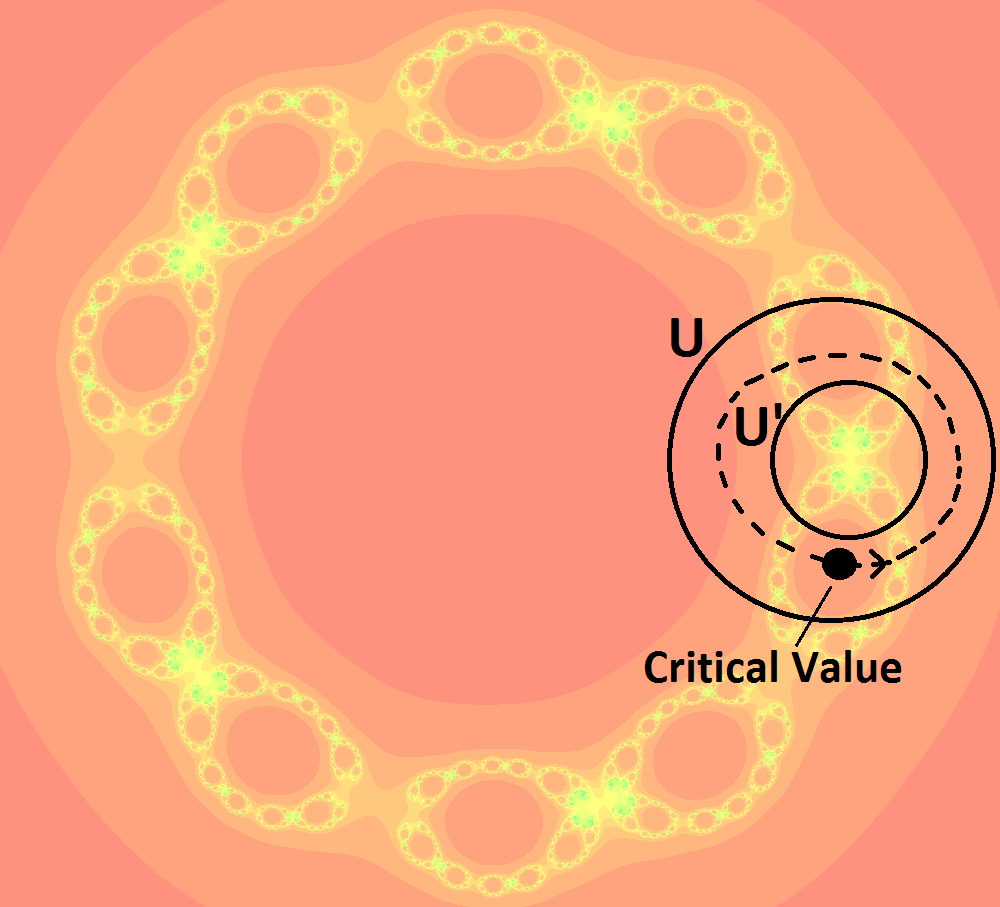}  
\caption{\label{fig:MandelwithW} Establishing the location of a baby Mandelbrot set requires tracking behavior in the dynamical plane as the parameter moves through the parameter plane.}
\end{figure}

Figure~\ref{fig:MandelwithW} illustrates a theoretical loop around a $W$ and corresponding loop around $U \setminus U'$.

This theorem is the method for proving where a given parameter plane contains homeomorphic copies of the Mandelbrot set.  Devaney used this method to show that the $a$-parameter plane of $R_{n,a,0}$ contains $n-1$ baby $\cM$'s for any $n \geq 3$. In fact, there appear to be many more than just $n-1$ copies, but these are the largest (and with the simplest critical orbit behavior) and so he refers to them as the  ``principal'' copies of the Mandelbrot set.

\subsection{Preliminaries Regarding $\Rnba$}

In this section, we collect some needed results on the family of functions $\Rnba(z)=z^n + \frac{a}{z^n} + b$ from \cite{BoydBrouwer1, BoydHoeppner1, BoydMitchell, BoydSchulz}.

\begin{notation} 
We use the following notation:
\begin{itemize}
    \item 
$\mathbb{D}(z_0,t) = \left\{z ~ \middle|~ \lvert z - z_0 \rvert < t\right\}$ is the disc of radius $t$ centered at $z_0$.
    \item 
$\mathbb{A}(t,s) = \left\{z ~ \middle|~ t < \lvert z \rvert < s \right\}$ is the annulus centered at the origin lying between the radii $t$ and $s$. 
    \item 
$\mathbb{A}(t,s)+ \gamma = \left\{z ~ \middle|~ t < \lvert z-\gamma \rvert < s \right\}$  is the annulus centered at $\gamma\in \CC$ and lying between the radii $t$ and $s$. 
\end{itemize}
\end{notation}

\begin{cor}\cite{BoydSchulz} \label{cor:BS1}
For any $b \in \CC$ and any $a \in \CC^{*}$, given any $\epsilon > 0$, there is an $N \geq 2$ such that for all $n \geq N$, we have that  the filled Julia set of $\Rnba$ is contained within an annulus near the unit circle: $K(\Rnba) \subset \mathbb{A}(1-\epsilon, 1+\epsilon)$.
\end{cor}

Of course, even in our subfamily $\Rna$ where $b$ is determined by $a$, this still applies, and in fact it will help us determine a neighborhood of the bifurcation locus for any choice of parameters, for $n$ sufficiently large.

 We also take advantage of a useful symmetry property of the family of functions $\Rnba$.
\begin{lem}\cite{BoydMitchell}
\label{lem:Involution}
Each map in the family of functions $\Rnba$ is symmetric under the involution $h_{a}(z)=\frac{a^{1/n}}{z}$.
\end{lem}

Additionally, we use some results from the first and last author's  \cite{BoydHoeppner1}.

\begin{prop} \cite{BoydHoeppner1}
\label{prop:JuliaAnnulus}
Let $b \in \CC$, $a \in \CC^*$, and $n \geq 3$. Let $s=\max\{4,|b|,|a|\}$ and $t=\dfrac{|a|^{1/n}}{s}$. Then 
 $K(\Rnba) \subset \mathbb{A}(t,s)$.
\end{prop}

We have pointed out that this family has $2n$ critical points and two critical values, $v_{\pm} = b \pm 2\sqrt{a}$.
In addition, $R$ maps half the critical points (counting multiplicity) to $v_+$ and the other half to $v_-$. 
\begin{lem} \cite{BoydHoeppner1}
The critical points of $\Rnba$ satisfy
\label{lem:critptsmapping}
\begin{equation}
\label{eqn:allcritpts}
a^{1/2n} = |a|^{1/2n}e^{i(\frac{\psi + 2k\pi}{2n})} 
\;\;\; \text{for} \;\;\; k=0,1,\ldots,2n-1,
\end{equation}
where $\psi= \Arg(a)$. Moreover, the $n$ critical points defined by $k$ even map to $v_+=b+2\sqrt{a}$ while the $n$ critical points with $k$ odd map to $v_-=b-2\sqrt{a}$.
\end{lem}
That is, $\crpt_{2j-1} = |a|^{1/2n}e^{i(\frac{\psi + 2(2j-1)\pi}{2n})}$
for $j=1,2,3,...,n,$ map to $v_-$ and 
$\crpt_{2j}=|a|^{1/2n}e^{i(\frac{\psi + 2(2j)\pi}{2n})} $ for
$j=0,1,2,...,n-1,$ map to $v_+.$ 

Note also that $\crpt_0 = |a|^{1/2n} e^{i\frac{\psi}{2n}}$ and 
 $\crpt_n = -|a|^{1/2n} e^{i\frac{\psi}{2n}}$.

\subsection*{The set $U_{a,k}$ and its image for any $\Rnba$.}

\begin{defn} \cite{BoydHoeppner1}
\label{defn:U_ak_prime}
For the family of functions $\Rnba$, let $\psi = \Arg(a)$ and for each $k=0,\ldots,2n-1,$ define 
$U_{a,k}'$ (which also depends on $n$),  
as the polar rectangle
 $$U'_{a,k} = \left\lbrace z \;\middle|\; \frac{|a|^{\frac{1}{n}}}{2} < |z| < 2, \; 
\left| \Arg(z) - \left( \frac{\psi+2k\pi}{2n} \right) \right| < \frac{\pi}{2n} 
\right\rbrace.$$
Also, set $U_a^+$ to be the image of $U'_{a,0},$ $U_a^-$ to be the image of $U'_{a,1}$, and use $U_{a,k}=R(U'_{a,k})$ when the parity of $k$ is not explicit:
$$U_{a}^+ = \Rnba (U'_{a,0}), \ \text{ and } \  U_{a}^- = \Rnba (U'_{a,1}),  \ \text{ and } \ U_{a,k}=R(U'_{a,k}).
$$
\end{defn}

So $U'_{a,k}$ is the polar rectangle with arguments centered at $(\psi + 2k\pi)/2n$ and then plus or minus $\pi/2n$.
Note $U'_{a,k}$ does not depend on $b$, but its image does.

\begin{defn} \cite{BoydHoeppner1}
\label{defn:ellipseEE}
 Let $\EE$ denote the ellipse parameterized by
$$
x(\theta) = \left( 2^n + {|a|}/{2}\right) \cos(\theta), \ \ \text{ and }
y(\theta) = \left( 2^n - {|a|}/{2}\right) \sin(\theta),
$$
then shifted so that it is centered at $b$ and rotated counter-clockwise by ${\psi}/{2}$.
\end{defn}

Note before the shift and rotation, the ellipse $\EE$ is centered at $0$, with its semi-major axis lying along the $x$-axis and of length $\left( 2^n + {|a|}/{2}\right)$, and its semi-minor axis lying along the $y$-axis and of length $\left( 2^n - {|a|}/{2}\right)$. 

By Lemma 6 of \cite{BoydMitchell}, the foci of $\EE$ are at the critical values $v_{\pm}.$ The images of the $U'_{a,k}$ sets were determined in \cite{BoydHoeppner1}:

\begin{prop} \cite{BoydHoeppner1}
\label{prop:Uellipse}
For each $k$  in $\{ 0,1, \ldots, 2n-1\}$, if $k$ is even,
$U_a^+ = \Rnba(U'_{a,k})$, 
and this set is one half of the ellipse $\EE$ including the minor axis and the critical value $v_+$. 
Moreover, $\Rnba$ maps 
$U'_{a,k}$ $2:1$ onto $U_a^+$ for each even $k$.

Similarly, 
$U_a^- = \Rnba(U'_{a,k})$ for each odd $k$
is the other half of the ellipse $\EE$, including the minor axis and the other critical value $v_-$, and $\Rnba$ maps 
 $U'_{a,k}$ two-to-one onto $U_a^-$ for each odd $k$.
\end{prop}

\subsection*{External Angle Assignments on a subset of $J(R_{n,a,b})$}

Since $R_{n,a,b}$ has two critical values, one of them may lie in a baby Julia set; some examples were constructed in \cite{BoydHoeppner1,BoydMitchell}. In \cite{BoydBrouwer1}, the first two authors provided a method for assigning external angles to the boundary of any baby Julia set that might exist, and to all of its preimages, in a way respecting the dynamics.

\begin{theorem}\cite{BoydBrouwer1} \label{thm: angle assignments}
    Let $R=R_{n,a,b}$ be a generalized McMullen map that is polynomial-like on a region $U'$ containing the critical value $v_+$, and that is conjugate on $K_+$, the filled Julia set of $R|_{U'}$, to a quadratic polynomial $P_c$ on its filled Julia set $K(P_c)$.
    Assume that $K(P_c)$ is connected and locally connected. Let $J_+ = \partial K_+$ and let $J_* = \cup_{m=0}^{\infty} R^{-m}(J_+)$.

    Then there exists a surjective relation $\Gamma:S^1\to J_*$ which assigns an angle in $[0,1)$ to each point in $J_*$ so that the angle assignments respect the dynamics to and from that point. In particular, if $J_{m,j}\in R^{-m}(J_+)$ is a preimage copy of $K_+$, then $\Gamma|^{J_{m,j}}$, $\Gamma$ restricted to co-domain $J_{m,j}$, is a surjective function $\Gamma|^{J_{m,j}}:S^1\to J_{m,j}$, and
    \begin{enumerate}
        \item If $m=0$, $R(\Gamma|^{J_+}(t))=\Gamma|^{J_+}(2t)$.
        \item If $m\geq1$ and $K_{m,j}\neq K_+$ is a component of $R^{-m}(K_+)$ which contains a critical point of $R$, then $R(\Gamma|^{J_{m,j}}(t))=\Gamma|^{R(J_{m,j})}(2t)$.
        \item If $m>1$ and $K_{m,j}$ is a component of $R^{-m}(K_+)$ which does not contain a critical point of $R$, then $R(\Gamma|^{J_{m,j}}(t))=\Gamma|^{R(J_{m,j})}(t)$.
    \end{enumerate}
\end{theorem}

The proof is inductive, first using the identification of $S^1$ with the interval $[0,1)$ and pulling those angles back to $J_+$ using the map $\phi$ that conjugates the dynamics of $R$ on $K_+$ to that of $P_c$ on $K(P_c)$, and using the dynamically significant external angle assignments on $K(P_c)$. Then those angle assignments can be pulled back through successive preimages of $J_+$ (whose union we called $J_*$) to assign angles on each (eventual) preimage copy of $J_+$ in a way that respects the dynamics.
Since $R$ is polynomial-like on $K_+$, $K_+$ must contain a critical value, so each of its direct preimages must contain a critical point. The local dynamics of $R$ are degree 2, so we see that the direct preimages of $J_+$ all map $2:1$ onto $J_+$, so the relation $\Gamma$ that assigns angles onto each preimage is conjugate to angle doubling on this domain. At each sequential set of preimages, angles are assigned based on whether the new preimages contain a critical point or not. If not, $\Gamma$ is conjugate to the identity map and angles are assigned as ``clones'' of their images. If a critical point is present, $\Gamma$ is again conjugate to angle doubling. 

\section{Dynamical Plane Results for $\Rna$}
\label{sec:dpr}

In this section, we describe how to build the subfamily $\Rna$, and provide some dynamical plane results for this subfamily. 

We begin with
$$
R_{n,a,b}(z) = z^n + \frac{a}{z^n} + b,
$$ 
then set $b$ so that $v_+ = b+2\sqrt{a}=a^{1/2n}$ holds for the canonical choice of root for $a^{1/2n},$ which we'll call the principle critical point. Specifically, if $\Arg(a)=\psi$, then we constrain $b$ by:
\begin{equation}
\label{eqn:vplus}
v_+ = b+2\sqrt{a}=|a|^{\frac{1}{2n}}e^{\frac{\psi}{2n}}.
\end{equation}
Solving for $b$, we see that if we set 
\begin{equation}
\label{defn:bn1a}
b_{n,0}(a) = |a|^{\frac{1}{2n}}e^{\frac{\psi}{2n}} - 2\sqrt{a},
\end{equation}
we are interested in the behavior of the family of functions
\begin{equation}
\label{eqn:Rna_definition_bfunction}
\Rna(z)= R_{n,a, b_{n,0}(a)}(z) = z^n + \frac{a}{z^n} + b_{n,0}(a),
\end{equation}
where $n \geq 3$ and $a \in \CC^*$. 
For simplicity, we refer to this parameterization as $\Rna$:
\begin{equation}
\label{eqn:Rna_definition_formula}
\Rna(z) =  z^n + \frac{a}{z^n} +|a|^{\frac{1}{2n}}e^{\frac{\psi}{2n}} - 2\sqrt{a}.
\end{equation}

We also note that since one critical orbit is always bounded,  we'll use $\Mn(\Rna)$ to refer to the set of parameters $a$ for which the free critical orbit is bounded. The boundary of this set is the bifurcation locus for this subfamily.

Now, the other critical value $v_- = b-2\sqrt{a}$ is free.
From Lemma~\ref{lem:critptsmapping}, we see that $\mapRna$ maps half the critical points (counting multiplicity) to $v_+$ and the other half to $v_-$. 


In this subfamily, the critical point with $k=0$ is equal to $v_+$. 

\smallskip

\noindent \textbf{Notation}.
Let $D_+$ denote the Fatou component containing $v_+$, and if $v_-$ is in the Fatou set, let $D_-$ denote the Fatou component containing $v_-$.

Since $v_+$ is a super-attracting fixed point here, $D_+$ is its immediate basin of attraction and $\mapRna(D_+) \subseteq D_+$.

\smallskip

\begin{figure}
\centering
\begin{subfigure}{0.5\textwidth}
  \centering
\includegraphics[width=.95\textwidth,keepaspectratio]{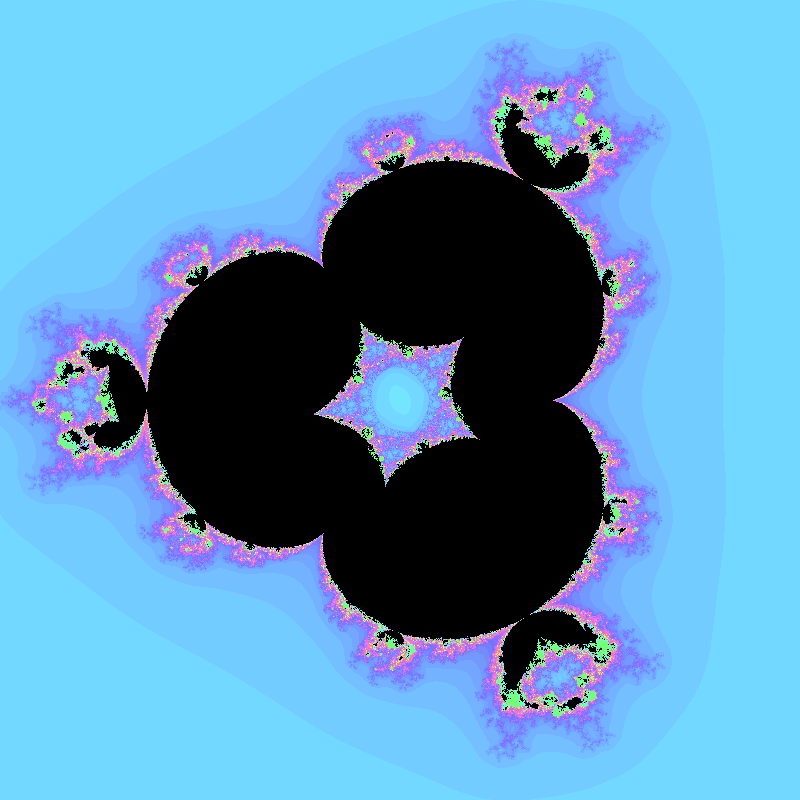}
  \caption{Parameter slice: $b$-plane for $R_{n,a,b}$ with $n=4, a=0.0001+0.0001i.$}
\end{subfigure}%
\begin{subfigure}{0.5\textwidth}
  \centering
\includegraphics[width=.95\textwidth,keepaspectratio]{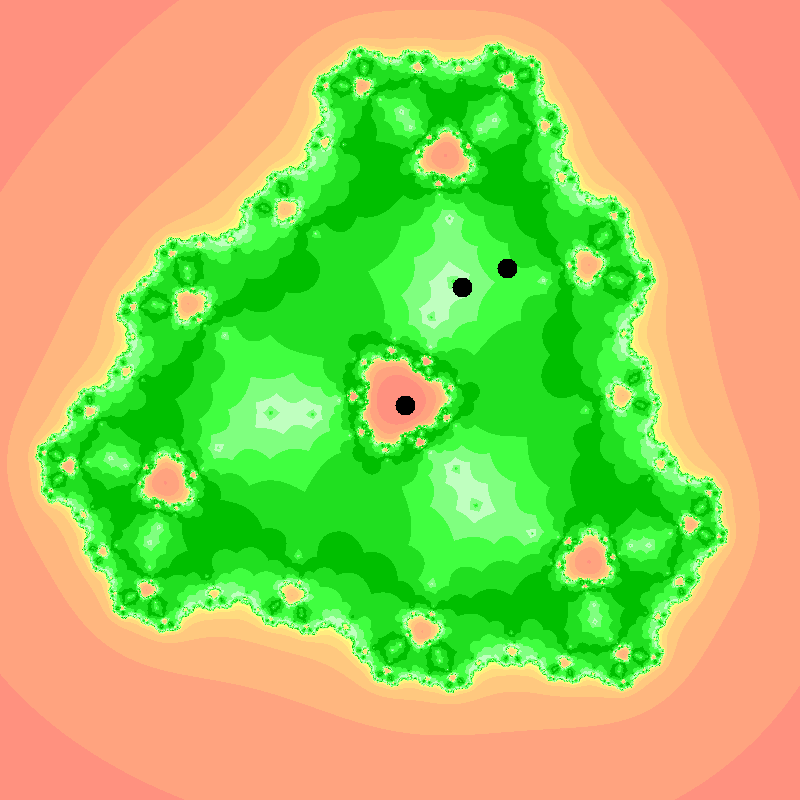}
  \caption{Dynamical plane for $n=3, a=0.001+0.001i, c=0.25+0.4i$, with critical values (and $0$) marked.}
\end{subfigure}%
\caption{Finding parameters for which both critical values lie in the same Fatou component.}
\label{fig:asmall}
\end{figure}

\noindent \textbf{Excluding a ``degenerate'' case.}
 For $n\geq 3$, if $a$ is small it can happen that $D_-=D_+$, that is, $v_-$ and $v_+$ are in the same (bounded) Fatou component. In that case, there is only one Fatou component, which the reader may imagine as a topological disk with an infinite tree of disks removed (the Fatou component containing $z=0$ and all of its preimages).   Xiao, Qiu, and Yin (\cite{xiaoqiu}) call such a $J$ a Cantor set of bubbles. 
In our subfamily with one critical point fixed, the hyperbolic component containing $a=0$ is a quite small teardrop-like shape aligned along the real axis, roughly in the bounds from $\Re(a)=-0.004$ to $\Re(a)=0.0025$ and $\Im(a)=-0.0025$ to $\Im(a)=0.0025.$ 
In the images in Figure~\ref{fig:FixedCritPointExamples1}, this component is not really visible as it lies just to the right of the ``deformed'' baby $\cM$-like set near the cusp of the cardioid shape of the bifurcation locus.  But in the general, parameter space maps with $D_-=D_+$ are quite easy to find. For example, on the left of Figure~\ref{fig:asmall} is the parameter slice which is the $b$-plane, holding $a$ small and constant, in this case for $n=4, a=0.0001+0.000i.$ The image appears like the parameter space for $z^4+b$, except with ``bites'' taken out of it due to the singular perturbation of $0$ mapping to $\infty$. Any point in the main central annulus-like black component  has both critical values in this common Fatou component. One such example is in this figure on the right, for $n=3, a=0.001+0.001i, c=0.25+0.4i$, although the situation for $n\geq3$ is similar. The critical values are marked, and the shading inside of the filled Julia set illustrates the rate of attraction to the attracting fixed point.  

In this article, we exclude this case of $D_-=D_+$, which geometrically is a small portion of our subfamily. For the rest of the parameter space, we have the following.

\begin{lem}
\label{lem:Ra-polynlikeD+}
Consider any $a$ in this parameter space such that $\mapRna=\Rna$ is hyperbolic, and $v_-$ and $v_+$ are not in the same Fatou component (so $D_-\neq \Dplus$).
Then there is a neighborhood $U=\cN(\overline{\Dplus})$ such that $\mapRna \colon U \cap \mapRna^{-1}(U) \to U$ is polynomial-like of degree two and conjugate to $p_0(z)=z^2$ on the set of points whose orbits never leave $U$. 
The conjugacy map is $\phi_a \colon U \to \cN(\overline{\mathbb{D}}) $, where $\phi_a$ is a quasi-conformal homeomorphism which is analytic on $\overline{\Dplus}$. Note $p_0 \circ \phi_a = \phi_a \circ \Rna$ and $\Dplus$ is a quasi-disk.
\end{lem}

\begin{proof}
Note by Proposition~\ref{prop:Uellipse}, $v_+$ is the only critical point in $\Dplus$.
Since $\mapRna$ is a hyperbolic rational map, there is a Riemannian metric for which $\mapRna$ is expanding on the boundary of $\Dplus$, which is a subset of the Julia set of $\mapRna$. Thus there is a neighborhood $U$ of $\overline{\Dplus}$ such that $\mapRna \colon U \cap \mapRna^{-1}(U) \to U$ is polynomial-like of degree two, and it must be conjugate to $p_0$ since its fixed point is superattracting. Note
$\phi_a$ is analytic on $\overline{\Dplus}$ by Theorem~\ref{thm:DH2}. 
\end{proof}

Thus, for most of the parameters in this parameter space, $\Dplus$ is not merely a simple topological disk, but a ``baby Julia set'', quasi-conformally homeomorphic to the unit disk the Julia set of $p_0,$ with the superattracting fixed point well ensconced inside the Fatou component, and, note $\partial \Dplus$ is a quasi-circle.

\begin{figure}
\centering
\begin{subfigure}{.5\textwidth}
  \centering
\includegraphics[width=.95\textwidth,keepaspectratio]{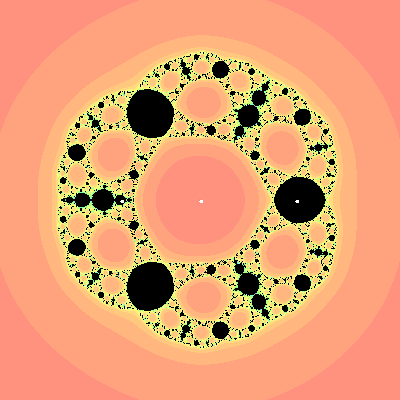}
  \caption{$a=0.094,$ chosen from a baby $\cM$}
  \label{Jfixedcrit_circlebas}
\end{subfigure}%
\begin{subfigure}{.5\textwidth}
  \centering
\includegraphics[width=.95\textwidth,keepaspectratio]{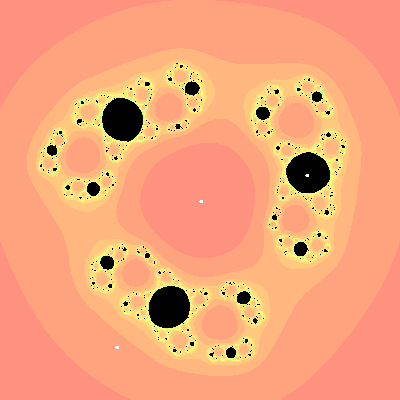}
  \caption{$a=0.02+0.2i$, from $\CC \setminus M_n(\Rna)$}
  \label{Jfixedcrit_outsideMn}
\end{subfigure}
\begin{subfigure}{.5\textwidth}
  \centering
\includegraphics[width=.95\textwidth,keepaspectratio]{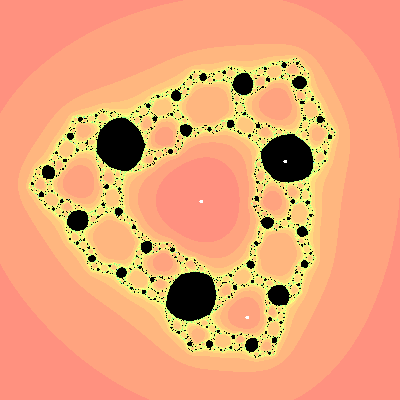}
  \caption{$a=-0.07+0.036i$ chosen from a Sierpinski hole}
  \label{Jfixedcrit_Sierpinskihole}
\end{subfigure}%
\begin{subfigure}{.5\textwidth}
  \centering
\includegraphics[width=.95\textwidth,keepaspectratio]{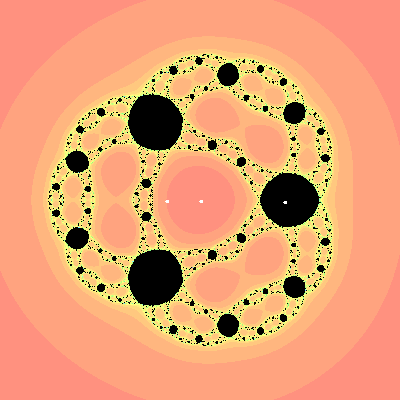}
  \caption{$a=0.043$ chosen from the central interior component of $\CC \setminus M_n(\Rna)$}
  \label{Jfixedcrit_CantorCirlces}
\end{subfigure}
\caption{Some Julia sets for $\Rna$, with $n=3$.
Critical values and $0$ are marked in white.}
\label{fig:FixedCritPointJulias}
\end{figure}

See Figure~\ref{fig:FixedCritPointJulias} for some examples of Julia sets for the case $n=3.$

\medskip

Next, we'll establish some bounds on the location of the filled in Julia set for maps in this subfamily. 

\begin{lem}
\label{lem:sizeofb}
For $a\neq 0$ set $b=b_{n,0}(a).$ Then 
$$
||a|^{1/2n} - 2|a|^{1/2}| \leq |b| \leq \sqrt{4|a| + |a|^{1/n} }.
$$    
\end{lem}

\begin{proof}
First observe that a short calculation shows that $|a|^{1/2n} \geq 2|a|^{1/2}$ when $|a| \leq (1/4)^\frac{n}{n-1}.$ Thus, we either have
\textit{Case 1:}  $|a| \leq (1/4)^\frac{n}{n-1}$ and the bounds will be $ 
0 \leq |a|^{1/2n} - 2|a|^{1/2} \leq |b| 
\leq 
  \sqrt{4(1/4)^\frac{n}{n-1} + (1/4)^\frac{1}{n-1} }
$, 
or \textit{Case 2:} $ (1/4)^\frac{n}{n-1} < |a|$ and we will have bounds
$ 0 <
2|a|^{1/2} - |a|^{1/2n} \leq
|b| \leq \sqrt{4|a| + |a|^{1/n} }.
$

Recall $b_{n,0}(a)= |a|^{\frac{1}{2n}}e^{\frac{\psi}{2n}} - 2\sqrt{a} = |a|^{\frac{1}{2n}}e^{\frac{\psi}{2n}} - 2|a|^\frac{1}{2}e^{\frac{\psi}{2}}.$
We consider the modulus of the sum of $|a|^{\frac{1}{2n}}e^{\frac{\psi}{2n}}$ and $- 2|a|^\frac{1}{2}e^{\frac{\psi}{2}}.$

Consider $a$ in the upper half plane, so as $\Arg(a)=\psi: 0 \to \pi,$ we have $\psi/2 = \Arg(2\sqrt{a}) : 0 \to \pi/2,$ and $\Arg(-2\sqrt{a}): -\pi \to -\pi/2$,  so $-2\sqrt{a}$ lies in quadrant 3. At the same time, $\Arg(a^{1/2n}): 0\to \pi/2n,$ so $a^{1/2n}$ lies in quadrant 1. Also, the argument of $-2\sqrt{a}$ changes faster than the argument of $a^{1/2n}.$
The minimum for $|b|$ occurs when $a$ is real and positive, so $a^{1/2n}$ is real and positive while $-2\sqrt{a}$ is real and negative. Then $\min |b| = | 2|a|^{1/2} - |a|^{1/2n}|$.
The maximum of $|b|$ occurs when $a$ is real and negative, in which case $-2\sqrt{a}$ lies on the negative imaginary axis and $a^{1/2n}$ is inside quadrant 1 (with argument $\pi/2n$). Then an upper bound can be found for $|b|$ by replacing $a^{1/2n}$ with $|a|^{1/2n}$ on the real axis, i.e. as $n$ may be large replace $\pi/2n$ with $0$, to see that  $\max|b|$ is bounded by the hypotenuse of the triangle with sides $|a|^{1/2n}$ and $2|a|^{1/2}$, that is,  $\sqrt{(2|a|^{1/2})^2 + (|a|^{1/2n})^2} 
= \sqrt{4|a| + |a|^{1/n} }.
$

The case of $a$ in the lower half plane is symmetrical. 
    
\end{proof}

Given this bound on $|b|$ for this subfamily, we can provide bounds on the location of the filled in Julia set in this subfamily.
 
\begin{lem}
\label{lem:K_annulus}
    For $n\geq 3$, $a\neq 0$ and $b=b_{n,0}(a),$ there exists an increasing sequence $\{q_n \}\subset (3.6, 4)$ and a decreasing sequence $\{\rho_n \} \subset (4, 4.4)$, where:
\begin{enumerate}
    
\item for $a$ such that $0 < |a| \leq q_n$, we have 
$$K(\Rna) \subset \mathbb{A} \left(\frac{|a|^{1/n}}{4}, 4 \right),$$ 
\item for $a$ such that $q_n \leq |a| \leq \rho_n,$ , we have $s =\sqrt{4|a| + |a|^{1/n} },$ hence 
$$K(\Rna)\subset \mathbb{A} \left(\frac{|a|^{1/n}}{\sqrt{4|a| + |a|^{1/n} }}, \sqrt{4|a| + |a|^{1/n} } \right),$$ 
\item for $a$ such that $\rho_n \leq |a|$, we have $s=|a|$,
hence 
$$K(\Rna)
\subset \mathbb{A} \left( \frac{|a|^{1/n}}{|a| }, |a| \right)
=\mathbb{A} \left( \left( \frac{1}{|a|} \right)^{\frac{n}{n-1}}, |a| \right).$$
\end{enumerate}

\end{lem}

\begin{proof}
Lemma~\ref{prop:JuliaAnnulus} states
for $b \in \CC$, $a \in \CC^*$, and $n \geq 3$, $s=\max\{4,|b|,|a|\}$ and $t=\dfrac{|a|^{1/n}}{s}$, we have 
 $K(\Rnba) \subset \mathbb{A}(t,s)$.

We combine this with Lemma~\ref{lem:sizeofb}. 
Since $|b| \leq \sqrt{4|a| + |a|^{1/n} }$, 
we get  $s\leq \max \{4, |a|, \sqrt{4|a| + |a|^{1/n} }\} $
and $t \geq |a|^{1/n} / (  \max \{4, |a|, \sqrt{4|a| + |a|^{1/n} }\})$.

Now we can compare the curves $y=4, y=x,$ and $y=\sqrt{4x+x^{1/n}}$ to determine $s$. There are three cases. 

\textit{Case 1:}
The curve $h_n(x)=\sqrt{4x+x^{1/n}}$ passes through the origin, is increasing with $x$ and is concave down. Thus, for $x>0$ small, $x < h_n(x) < 4$. Then for each $n\geq 3$ there is a $q_n$ such that
$h_n(q_n) = 4,$ and  as $n$ grows from $3$ to $\infty$, $q_n $ grows from approximately $3.6163$ toward $4$. 
Note that $(1/4)^{n/(n-1)} < q_n$ (in fact, as $n$ grows from $3$ toward $\infty$, $(1/4)^\frac{n}{n-1}$ grows from $1/8$ toward $1/4$).
Thus, in the interval 
$0 < |a| < q_n, $ we have $|a| < h_n(|a|) < 4,$ and hence $s= \max \{4, |a|, \sqrt{4|a| + |a|^{1/n} }\})=4$.

\textit{Case 2:}
Next, for each $n\geq 3$ there is a $\rho_n$ such that 
$h_n(\rho_n) = \rho_n$.
As $n$ grows from $3$ toward $\infty$, $\rho_n$ decreases from about $r_3\approx 4.37390$ toward $\lim_{n\to\infty} \rho_n = 4$. 
Hence, $\rho_n > 4$ for all $n$.
Thus, in the interval
$q_n < |a| < \rho_n,$ (a little less than $4$ to a little more than $4$), we have $s =\sqrt{4|a| + |a|^{1/n} }\}. $ 

\textit{Case 3:}
Then finally when $4 < |a|$, we have $s=|a|.$
Then $t=|a|^{1/n} / |a| = (1/|a|)^{n/(n-1)}.$

\end{proof}

We can also get some slightly different bounds on the location of the Julia set that will be useful for parameters in or near the boundedness locus. 

\begin{lem}
    \label{lem:Kbounds_a_in_unitdisk}
If $|a|<1$, for all $n \geq 3$, we have $K(\Rna) \subset  \mathbb{A} \left(\frac{|a|^{1/n}}{2}, 2 \right) \subset \mathbb{A} \left(\frac{1}{2}, 2 \right)$.
\end{lem}

\begin{proof}
    The proof of Lemma 2 in \cite{BoydSchulz} says that if $N$ satisfies $(1+\ep)^N > 3\max(1, |a|, |b|)$, then for all $n \geq N$, the escape radius is $(1+\ep)$. Set $\ep=1$, so we need $2^3=8 > 3\max (1, |a|, |b|),$ but $|a|<1$ so we need $8> \max (3, 3|b|),$ hence we need $8>\max(3|b|).$ In this subfamily, $b=|a|^{1/2n}e^{i\psi/2n}-2\sqrt{a}.$
    
So for $|a|<1$, consider the case that $\Im(a)\geq 0$, so $a$ is in the upper half plane. Then first, the principal root $a^{1/2n}$ lives in the wedge in quadrant $1$ with modulus from $0$ to $1$ and with argument going from $0$ to $\pi/(2n) \leq \pi/6$ as the argument of $a$ goes from $0$ to $\pi$.   Second, $-2\sqrt{a}$ lives in the quadrant $3$ (both real and imaginary parts non-positive) intersect the disk of radius $2$, with argument going from $-\pi$ ``up'' to $-3\pi/2$ as the argument of $a$ goes from $0$ to $\pi$. Now considering the sum 
$b=a^{1/2n} + (-2\sqrt{a})$, 
the maximum is bounded above by the sum when $a=-1=e^{i\pi}$, 
 so $-2\sqrt{a}=-2i$, and $a^{1/2n}=(e^{i\pi})^{1/2n}$ is at the upper right corner of the wedge. Thus the maximum over all $n$ of $(-1)^{1/2n}$ is bounded above by taking $n \to \infty$ so we have $(e^{i\pi})^{1/2n}=1$. Thus 
$|b| < |1-2i|=\sqrt{5}$. Fortunately, $\sqrt{5} < 2.24 < 8/3 = 2 \frac{2}{3},$ so $3|b| < 3\sqrt{5} < 8=2^3.$

Of course, the case that $a$ has negative imaginary part is symmetrical. Hence $|a|<1$ implies for all $n\geq 3$, $K(\Rna)\subset \DD(0,2).$ But then the involution property of the family (Lemma~\ref{lem:Involution}) implies that $K(\Rna)\subset \mathbb{A}(\frac{|a|^{1/n}}{2}, 2).$

\end{proof}

Next, as our free critical value $v_-$ is significant, the following is also useful. 

\begin{lem}
\label{lem:outerboundaryW} If $a$ satisfies $|v_-|\leq 2$, and $n\geq 3$, then $\DD(0,2)\subset \EE$, and hence $U'_{a,k} \subset \EE$ for all $k$. 
\end{lem}

\begin{proof}
As specified in Lemma~\ref{prop:Uellipse}, the length of the semi-major axis of $\EE$ is given by $2^n + \frac{|a|}{2^n}$ and the length of the semi-minor axis of $\EE$ is given by $2^n - \frac{|a|}{2^n}$. Recall  $\EE$ is centered at $b$ and is rotated by a degree of ${\psi}/{2}$, and the two foci of $\EE$ occur at $b \pm 2\sqrt{a}$, the two critical values of $\Rnba$. Hence, the distance from the center of the ellipse to either focus is $2\sqrt{|a|}$.

Now that we know the locations of the two foci of $\EE$ we can prove that $\overline{\mathbb{D}(0,2)} \subset \EE$ by examining the sum $|z - v_+| + |z - v_-|$ for $|z| \leq 2$. The sum of the distances from any point on $\EE$ to each foci is $2 \left( 2^n + \frac{|a|}{2^n} \right)$. If we can show  $|z - v_-| + |z - v_+| < 2 \left( 2^n + \frac{|a|}{2^n} \right)$ for all allowable choices of $a$, this would give us that $\overline{\mathbb{D}(0,2)} \subset \EE$. 

But we know $v_+$ is a fixed critical point, so $|v_+| = |a|^{1/2n},$ and we assumed $|v_-|\leq 2$. So,
\begin{align*}
& |z - v_+| + |z - v_-| 
\leq 2|z| + |v_+| + |v_-| 
\leq 4 + |a|^{1/2n} + 2 = 6 + |a|^{1/2n}.
\end{align*}
So we need $6 + |a|^{1/2n} < 2 \left( 2^n + \frac{|a|}{2^n} \right)$.
For ease of notation let $x=|a|^{1/2n}$, so $x^{2n}=|a|$ and $x>0$. So we need $6 + x < 2 (2^n + x^{2n}/2^n) = ,$ equivalently:
$
x^{2n} - 2^{n-1}x + 2^n(2^n-3) > 0,
$
for all $x>0$ and $n\geq 3$. 
Let $g_n(x) = x^{2n} - 2^{n-1}x + 2^n(2^n-3).$ 
Examination of $g_n(x)$ yields $g_n'(x) = 2n x^{2n-1} -2^{n-1},$ and
$g''(x)=2n(2n-1)x^{2n-2}$ always positive. So $g_n$ is concave up, and solving for $g_n'(x)=0$ we see the unique minimum of $g_n$ occurs at the value 
$x_n = \left(\frac{2^{n-2}}{n}\right)^{\frac{1}{2n-1}} >0,$ where (using some calculus) we get $x_n \uparrow \sqrt{2}$ as $n\to\infty$. Moreover, the derivative changes slowly near $x=0$, and $g_n(0)=2^n(2^n-3)$ is large. As a result, the function is not only positive but extremely so. In particular, 
as $g_n$ is concave up, the steepest slope of $g_n$ on $[0,x_n]$ occurs at $0$, and is $g'_n(0)=-2^{n-1}$. Thus, $g_n(x_n) \geq g_n(0) + g_n'(0)x_n = 2^n(2^n-3) - 2^{n-1} x_n. $ Using $x_n \leq \sqrt{2}$ and $n\geq 3,$ we get $g_n(x_n) \geq 2^n(2^n-3) - 2^{n-1}2^{1/2} \geq 2^n ( 2^n - (3+\sqrt{2}/2)) > 2^n (2^n - 4) \geq 32 >  0$. 

Thus $\overline{\mathbb{D}(0,2)} \subset \EE$ and therefore $U'_{a,k} \subset \EE$ (for such $a$).

\end{proof}
 

\section{Locating the bifurcation locus of $\Rna$}
\label{sec:spines}

Now we look more globally for a region in the parameter space in which the bifurcation locus $M(\Rna)$ lies. 

First, for $\Rna$ we are able to find a spine by examining the values of $a$ for which the bounded critical orbits lie on the unit circle, motivated by Corollary~\ref{cor:BS1}.

\begin{figure}
\centering
\includegraphics[scale=.2]{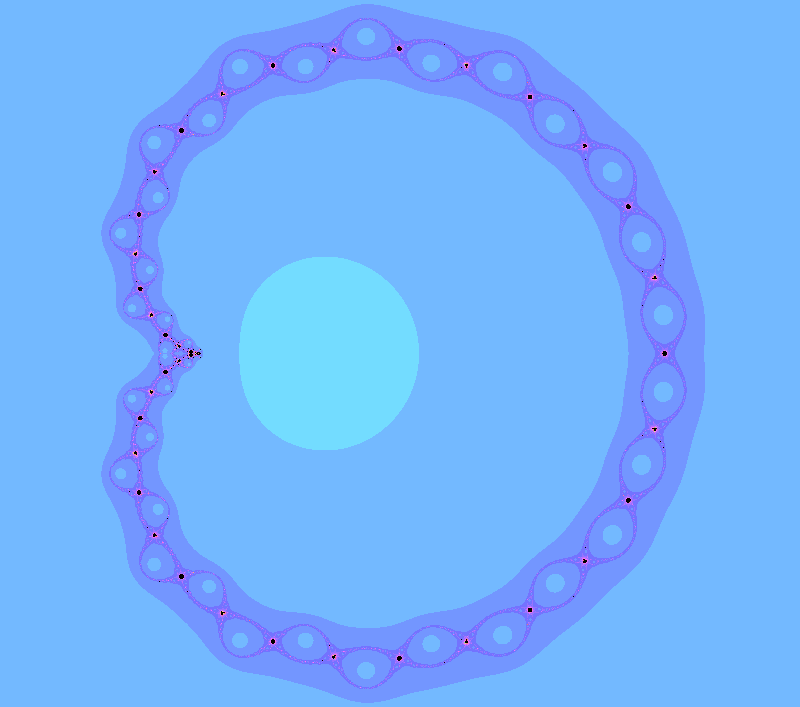}
\caption{Parameter plane for $\Rna$ for $n=20$}
\label{fig:fixed_crit_n=20_aplane}
\end{figure}

\begin{lem}\label{lem:spinecritfixed}
For $b=b_{n,0}(a)$, the set of all $a \in \CC$ for which at least one of the equations $|b \pm 2\sqrt{a}| = 1$ holds is given by the set of solutions:
\begin{equation*}
\cS_n
= \left\{ a=\frac{1}{16}(|a|^{1/2n} e^{i\frac{\psi}{2n}} + e^{i\theta})^2 \ \  | \ \  0\leq \theta \leq 2\pi
\right\}.
\end{equation*}
As $n$ gets large, this approaches the cardioid
\begin{equation*}
\cS_{\infty} = \left\{
a=\frac{1}{16}(1 + e^{i\theta})^2 \ \  | \ \  0\leq \theta \leq 2\pi
\right\},
\end{equation*}
which has a cusp at $0$, and real max at $1/4$.
\end{lem}
\begin{proof}
We first write $|b \pm 2\sqrt{a}| = 1$ as $b \pm 2\sqrt{a} = -e^{i\theta}$ for $0 \leq \theta \leq 2\pi$. Substituting $b=|a|^{1/2n} e^{i\frac{\psi}{2n}}-2\sqrt{a}$ gives us both:
(i): $|a|^{1/2n} e^{i\frac{\psi}{2n}} =  -e^{i\theta},$
and
(ii): $|a|^{1/2n} e^{i\frac{\psi}{2n}} - 4\sqrt{a} =  -e^{i\theta}.$
The first equation represents the critical value $b+2\sqrt{a}$, which we know is fixed and so will only lie near the unit circle when $|a|=1$.
As for the second equation, we get the form that we are looking for by isolating the square root term and then solving for that $a$.
\begin{align*}
& |a|^{1/2n} e^{i\frac{\psi}{2n}} - 4\sqrt{a} =  -e^{i\theta} \ \Rightarrow \ 
4\sqrt{a} = a|a|^{1/2n} e^{i\frac{\psi}{2n}} + e^{i\theta} 
\\ & 
\Rightarrow
\sqrt{a} = \frac{1}{4}|a|^{1/2n} e^{i\frac{\psi}{2n}} + e^{i\theta})
\  \Rightarrow  \
a = \frac{1}{16}(|a|^{1/2n} e^{i\frac{\psi}{2n}} + e^{i\theta})^2,
\end{align*}
and note the prinicipal root of $a^{1/2n}$ approaches $1$ as $n\to \infty$. 
\end{proof}
It may not be obvious from Figure~\ref{fig:FixedCritPointExamples1} that this spine becomes more like a cardioid as $n$ gets larger. If we trace the large black shapes within the necklace structure we can see a curve which looks at least somewhat like a cardioid as early as $n=6$. Figure \ref{fig:fixed_crit_n=20_aplane} shows an example of a large value of $n$ where we very clearly see a cardioid shape.

We next show the bifurcation locus lies within an annulus. 
In the limit of the above cardioid as $n\to\infty$, the cardioid has cusp at $a(\pi)=0$, and max on the $x$-axis at $a(0)=1/4$, so is contained in an annulus centered at $1/8$, with a small inner radius  $<1/8$, and outer radius about $1/7$. We show below is that the boundedness locus is contained in the annulus $\mathbb{A}(l,u)+1/8$, with center $1/8$, outer radius any $u > 3/8$ and inner radius any $l\leq 1/32,$ for $n$ sufficiently large.

\begin{lem}
\label{lem:M_in_annulus}
    Let $l \leq 1/32$ and $3/8 < u$. Then $\exists N\geq 3$ such that  $\forall n \geq N, M_n(\Rna) \subset \mathbb{A}(l,u)+{1}/{8},$  i.e., the set of parameters with the orbit of $v_-$ bounded 
is contained in the annulus centered at $1/8$ with inner radius $l$ and outer radius $u$.
\end{lem}

\begin{proof}
Note $\mathbb{A}(l,u)+{1}/{8} = \{ z \colon \ l < | z - 1/8| < u \} = \DD(1/8, u) \setminus \DD(1/8,l).$

Now, in this family we always have $|v_+| = |a|^{1/2n},$ so $v_+$ has a bounded orbit, so we'll focus on $v_-$, and $v_- = 
 |a|^{1/2n}e^{i\frac{\psi}{2n}} - 4\sqrt{a}$.

\textit{Outer}: Let $u=\frac{1}{8}+ \frac{(2+\delta)^2}{16} $ for some $0 < \delta< 2$. Suppose $a \in \CC \setminus \DD(1/8, u),$ so $|a-1/8| \geq u,$ hence $|a| \geq u-1/8$.
So,
$|v_-| = 
| |a|^{1/2n}e^{i\frac{\psi}{2n}} - 4\sqrt{a}| =
| 4\sqrt{a} - |a|^{1/2n}e^{i\frac{\psi}{2n}}| \geq
4 |a|^{1/2} - |a|^{1/2n}  = 
|a|^{1/2n}(4|a|^{\frac{n-1}{2n}} -1)
\geq 
|u-1/8|^{1/2n}(4|u-1/8|^{\frac{n-1}{2n}} -1).
$ Now, $u-1/8 = (2+\delta)^2/16 \in (1/4,1)$ for $\delta< 2$.
Next, $h(x,n)=x^{1/2n}(4x^{(n-1)/2n}-1)$ for fixed $0< x < 1$  is a decreasing function of $n$ for $n\geq 0$,
but $\lim_{n\to\infty} h(x,n) = 4\sqrt{x}-1.$
Hence, $h(u-1/8,n) > 4\sqrt{u-1/8}-1 =
4\sqrt{\frac{(2+\delta)^2}{16}}-1
= 2+\delta-1 = 1+\delta > 1.$
Thus, $|v_-| > 1 +\delta.$

Though we are more interested in a tighter bound for our outer annulus, we note a modification of our argument also yields $|v_-|> 1$ if $\delta\geq 2$. Indeed, one can check that $h(x,n)=x^{1/2n}(4x^{(n-1)/2n}-1)$ for fixed $x\geq 1$  is an increasing function of $n$ for $n\geq 2$, 
so $h(x,n)\geq h(x,3)=x^{1/6}(4x^{1/3}-1).$
But also, $h(x,3)$ is a  non-decreasing function of $x$ for all $x\geq 1$ (increasing if $x>1$, flat if $x=1$) so $h(x,3) \geq h(1,3)=3$ for all $x\geq 1$. 
Hence from the above, assuming $u-1/8\geq 1$ (which it is if $\delta \geq 2$) we have 
$|v_-| \geq 
|u-1/8|^{1/6}(4|u-1/8|^{\frac{1}{3}} -1) \geq 3
$.

So, for $a \in \CC \setminus \DD(1/8, u),$ as long as $u > 3/8$, we have 
$|v_-| \geq 1 + \delta'$ for some $\delta'>0$.  
 This means that by Corollary~\ref{cor:BS1}, for any choice of $a \in \CC \setminus \DD(1/8, u),$ and $u> 3/8$, there exists an $N_1 \geq 3$ s.t.\ for all $n\geq N_1$, we have 
 $v_- \notin K(\Rna)$, that is, $v_-$ is not in the filled Julia set, at least for $n$ large, hence $a\notin M_n(\Rna)$. 
 
\textit{Inner}: Suppose $a \in \DD(1/8,l),$ so $|a-1/8| < l.$ 

Say $l\leq 1/32.$ Then easily $\Arg(a) \in (-\pi/10,\pi/10)$, (in fact $\pi/12$ works) and $|a| \in (1/8-l, 1/8+l) \subseteq (1/8-1/32, 1/8+1/32) = (3/32, 5/32)$.  

First, $\sqrt{a}$ has argument in $(-\frac{\pi}{20}, \frac{\pi}{20})$ and modulus in $(\frac{\sqrt{3}}{4\sqrt{2}}, \frac{\sqrt{5}}{4\sqrt{2}})$, 
so $-4\sqrt{a}$ has argument in 
$\pi \pm \frac{\pi}{20},$ or $(-\pi, -\frac{19\pi}{20}) \cup (\frac{19\pi}{20}, \pi]$, 
and modulus in $(\sqrt{\frac{3}{2}}, \sqrt{\frac{5}{2}})$.

Second, since $n\geq 3$, $a^{1/2n}$ 
has argument in 
$(-\frac{\pi}{20n},\frac{\pi}{20n}) \subseteq (-\frac{\pi}{60},\frac{\pi}{60}),$ and modulus in $((\frac{3}{32})^{1/2n}, (\frac{5}{32})^{1/2n}) \subseteq ((\frac{3}{32})^{1/6}, (\frac{5}{32})^{1/6})$. 

Finally, adding these up we want to find an upper bound for $|v_-|$, call it $L(n)$. 
Since we are summing two polar rectangles with a larger one (of larger modulus and angle width) centered on the negative real axis and a smaller one centered on the positive real axis, the maximum modulus occurs when
$-4\sqrt{a}$ is in the upper left corner of its polar rectangle, and $a^{1/2n}$ is in the upper right of its polar rectangle (or the symmetrical choice). See Figure~\ref{fig:polar_rectangles_1}.
\begin{figure}
\includegraphics[scale=.45]{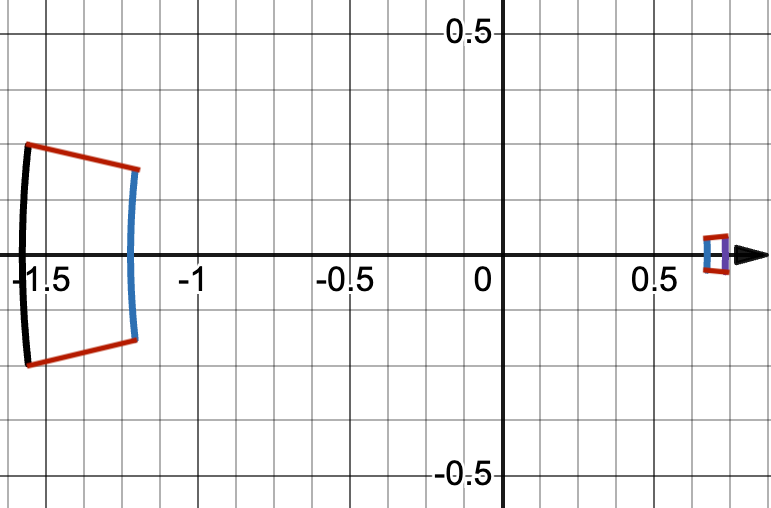}
\caption{\label{fig:polar_rectangles_1}Polar rectangles bounding $a^{1/2n}$ (right polar rectangle) and $-4\sqrt{a}$ (left polar rectangle) when $n=3$ and $|a-1/8| = 1/32$.
}
\end{figure}

This is
for $n=3$: $L(3)= | (\sqrt{5/2})e^{i(19\pi/20)} + (3/32)^{1/6}e^{i(\pi/60)}| < 0.95$, or in general at
$$L(n) = \left| \sqrt{\frac{5}{2}}\ e^{i\frac{19\pi}{20}} + \left(\frac{3}{32}\right)^{\frac{1}{2n}}e^{i\frac{\pi}{20n}}\right|,$$ which is a decreasing function of $n$, less than $1$ when $n > 3$, with $L(4)< 0.87$, 
and in the limit as $n\to \infty$ we get $L(n)$ bounded above by $0.614.$
A few more bounds for $L(n)$ are given in Table~\ref{table:Ln}.

\begin{table}
\begin{tabular}{|c|c|c|c|c|c|c|c|c|c|c|c|}
    \hline
    $n$ & $3$ &  $4$ & $5$ & $6$ & $7$ & $10$ & $15$ & $25$ & $50$ & $100$ & $\infty$ \\ 
    \hline
       $L(n) < $  & $0.95$
       & 0.87 & $0.82$ & 0.8 & 0.77 & 0.73 & 0.7 & 0.66 & 0.64 & 0.63 & 0.62
       \\
    \hline
\end{tabular}
\caption{\label{table:Ln} Upper bounds for  $|v_-|$ for some values of $n$, if $|a - 1/8 | < 1/32$.}
\end{table}

Now considering Corollary~\ref{cor:BS1}, 
we know as $n$ increases, the filled Julia set tends to the unit circle, hence for $n> N_2$ sufficiently large, $L(n)$ is sufficiently less than $1$ that $v_-$ with modulus less than $L(n)$ is guaranteed not to have a bounded orbit. 

Thus taking $N = \max\{N_1,N_2\},$ we have for all $n\geq N, M_n(\Rna) \subset \mathbb{A}(l,u)+1/8.$
\end{proof}

The above quickly yields: 

\begin{corollary}
\label{cor:M_in_disk}
For any $\ep>0$, there is an $N\geq 3$ s.t.\ for all $n>N$, we have $M_n(\Rna) \subset \DD(0,1/2+\ep).$
\end{corollary}

In fact, computer images suggest the boundedness locus is contained in $\DD(0,1/3)$ for all $n\geq 3$. 

Now, using these estimates we can describe when $U' \subset \EE$.

\begin{corollary}
    \label{cor:Uprime_in_E}
For all $n$ sufficiently large, if $v_-\in K$, then $U'_{a,k} \subset \cE$.
\end{corollary}

\begin{proof}
Using $\ep=1/2$ in Corollary~\ref{cor:M_in_disk}, we have $M_n(\Rna)\subset \DD(0,1)$ for $n$ sufficiently large (computer images suggest $n=3$ suffices). Thus, for $a\in M_n(\Rna)$ we know $|a|<1$. So, by Lemma~\ref{lem:Kbounds_a_in_unitdisk}, we have $K(\Rna)\subset \mathbb{A}(|a|^{1/n}/2, 2).$
Note also then that for $a\in M_n(\Rna),$ (still for $n$ suff.\ large), we have if $v_-\in K $, then since $ K \subset \DD(0,2)$ we get $|v_-|< 2$. Thus Lemma~\ref{lem:outerboundaryW} applies, so $U'_{a,k} \subset \cE$ for all $k$, if $v_-\in K$ and  $n$ is so large that 
$M_n(\Rna)\subset \DD(0,1)$.
    \end{proof}

\subsection*{Future work on the location of the boundedness locus}
For the subfamily $\mapRna$, we have now identified a spine for the boundedness locus (Lemma~\ref{lem:spinecritfixed}), and showed the boundedness locus was contained in an annulus for sufficiently large $n$ (Lemma~\ref{lem:M_in_annulus}). 

One way to improve on Lemma~\ref{lem:M_in_annulus} would be to
establish the following: 
\begin{conjecture}
    $\forall \ep>0$, $\exists N\geq 3$ s.t.\ $\forall n \geq N, M_n(\Rna) \subset \cN_{\ep}(\spine)$ and/or  $M_n(\Rna) \subset \cN_{\ep}(\cS_\infty)$; that is, the boundedness locus
is contained in an $\ep$-neighborhood of the spine. 
\end{conjecture}

The outline of a proof could be somewhat similar to the proof of Theorem 2 given in \cite{BoydHoeppner1} (or to the proof of a simpler case, for a lemma in \cite{BoydSchulz}): show that inside the annulus guaranteed by Lemma~\ref{lem:M_in_annulus}, but outside of $\cN_{\ep}(\spine),$ $|v_-|$ is bounded away from $1$, then apply Corollary~\ref{cor:BS1}. The fact that the spine depends on $n$ is an added complication that does not appear in  \cite{BoydHoeppner1}, but computer images suggest $M_n(\Rna) \to \cS_{\infty}$ as $n\to\infty$.  

\section{Centers of Principal Components in the $\Rna$ Parameter Space}
\label{sec:centers}

Baby $\cM$'s naturally tend to appear around values of the parameter $a$ for which a free critical value happens to be fixed. For $\Rna$, since $v_+$ is always fixed,  first we need to locate values of $a$ for which the second critical value, $v_-$, is fixed as well. 

Fixing $v_- = b-2\sqrt{a}$ means we have $b-2\sqrt{a} = a^{1/2n}$ for one of the $2n^{\text{th}}$ roots of $a$. The critical value $b+2\sqrt{a}$ is already set equal to the canonical root $|a|^{\frac{1}{2n}}e^{\frac{\psi}{2n}}$.

This gives us one trivial solution: If $a=0$, then $b-2\sqrt{a} = b+2\sqrt{a} = a^{1/2n} = 0$. For the rest, we need $v_-$ to be equal to a non-canonical choice of root. 

\begin{lem}
\label{lem:akdefn}
The set of $a$-values for which the critical value $v_-=b-2\sqrt{a}$ is equal to one of the non-canonical roots of $a^{1/2n}$ is 
$$
\left\{ a = a_k = \left( \frac{1-e^{i\frac{k\pi}{n}}}{4} \right)^{\frac{2n}{n-1}} \colon   k=1,2,3,...,2n-1.\right\} 
$$

At $a=a_k$, $v_- = w_k = |a|^{1/2n}e^{i(\frac{\psi + 2k\pi}{2n})}$.
\end{lem}
\begin{proof}
Setting $v_- = b-2\sqrt{a}$ equal a root of $a^{1/2n}$ and substituting in our choice of $b= |a|^{1/2n} e^{i\frac{\psi}{2n}} - 2\sqrt{a}$ gives us: 
$$
 a^{1/2n} = b- 2\sqrt{a} =  |a|^{\frac{1}{2n}}e^{\frac{\psi}{2n}} - 2\sqrt{a} - 2\sqrt{a} = |a|^{\frac{1}{2n}}e^{\frac{\psi}{2n}} - 4\sqrt{a}  .$$
 In order to avoid a trivial equation, the $a^{1/2n}$ root on the left-hand side of the equation must \textbf{not} be the same canonical root. It must be chosen from one of the $2n-1$ other roots. Thus we really have $2n-1$ equations to consider:
$$|a|^{\frac{1}{2n}}e^{\frac{\psi}{2n}} 
- 4\sqrt{a} = 
|a|^{1/2n}e^{i(\frac{\psi + 2k\pi}{2n})} 
\;\;\; \text{for} \;\;\; k=1,2,3,...,2n-1.$$

We next isolate the square root, then solve for $a = |a| e^{i\psi}$:
\begin{align*}
& 4\sqrt{a} = |a|^{\frac{1}{2n}}e^{\frac{\psi}{2n}}  - 
|a|^{1/2n}e^{i(\frac{\psi + 2k\pi}{2n})} 
 \ \Rightarrow \ 
 \sqrt{a} = |a|^{\frac{1}{2n}}e^{i \frac{\psi}{2n}}\left( 
(1 - e^{i\frac{2k\pi}{2n}})/{4} \right) \ 
 \\ &  \Rightarrow 
a^{\frac{n-1}{2n}} = (1 - e^{i\frac{k\pi}{n}})/{4}
 \ \Rightarrow \ 
a = \left( (1 - e^{i\frac{k\pi}{n}})/{4} \right)^{\frac{2n}{n-1}}.
\end{align*}
\end{proof}

Around these $2n-1$ points, we see two patterns emerge. Parameters defined by the above with $k$ even lie inside the disks we see in the bifurcation locus. The $a_k$ for $k$ odd lie within observed baby $\cM$'s, except that $a_1$ and $a_{2n-1}$ lie in the same merged/deformed baby $\cM$ near $a=0$, which is not a true baby $\cM$. 
To begin to describe why, we examine the images of the critical value $v_-$ under $\Rna(z)$ at the $2n-1$ $a_k$-values listed above.

\begin{prop}\label{prop:c-patterns}
Consider $a = a_k = \left(  ({1-e^{i\frac{k\pi}{n}}})/{4} \right)^{\frac{2n}{n-1}} $ for some 
$k\in \{1,2,3,...,2n-1\},$ 
so $v_-$ is one of the non-principle critical points $w_k$.

If $k$ is odd then $\Rna(v_-)=v_-$, i.e., the critical value $v_-$ is a fixed point of $\Rna$. If $k$ is even, then $\Rna(v_-)=v_+$.
\end{prop}
\begin{proof}
For ease of notation, set $\alpha_k =(1-e^{i\frac{k\pi}{n}})/{4}$.
 We want 
$\Rna(v_-) = b\pm 2\sqrt{a},  $ equivalently,
$\Rna(v_-) - b  =  2\sqrt{a}$. That is, we want
$$
(b-2\sqrt{a})^n + \frac{a}{(b-2\sqrt{a})^n}  = \pm 2\sqrt{a}.$$
Using $b=a^{1/2n} - 2\sqrt{a}$ and $a = (\alpha_k)^{\frac{2n}{n-1}}$ in the left hand side above yields
$$\left( \alpha_k^{\frac{1}{n-1}} - 4\alpha_k^{\frac{n}{n-1}}\right)^n + \frac{ \alpha_k^{\frac{2n}{n-1}}}{\left( \alpha_k^{\frac{1}{n-1}} - 4\alpha_k^{\frac{n}{n-1}}\right)^n}.$$
Next we factor $\alpha_k^{\frac{n}{n-1}}$ out of both of these terms:
$$\left( \frac{\left( \alpha_k^{\frac{1}{n-1}} - 4\alpha_k^{\frac{n}{n-1}}\right)^n}{\alpha_k^{\frac{n}{n-1}}} + \frac{\alpha_k^{\frac{n}{n-1}}}{\left( \alpha_k^{\frac{1}{n-1}} - 4\alpha_k^{\frac{n}{n-1}}\right)^n}\right)\alpha_k^{\frac{n}{n-1}}.$$
Hence:

$$\left( \left(\frac{
\alpha_k^{\frac{1}{n-1}} - 4\alpha_k^{\frac{n}{n-1}}
}{
\alpha_k^{\frac{1}{n-1}}
}\right)^n 
+ 
\left(\frac{
\alpha_k^{\frac{1}{n-1}}
}{
\alpha_k^{\frac{1}{n-1}} - 4 \alpha_k^{\frac{n}{n-1}}
}\right)^n \right)
\alpha_k^{\frac{n}{n-1}}.$$

After some cancellation we have:
$$\left( \left( 1 - 4\alpha_k^{\frac{n-1}{n-1}}\right)^n + \left(\frac{1}{1 - 4\alpha_k^{\frac{n-1}{n-1}}}\right)^n \right)
\alpha_k^{\frac{n}{n-1}}.$$
Simplifying further leaves us with
$$\left( e^{ik\pi} + e^{-ik\pi} \right) \alpha_k^{\frac{n}{n-1}}.$$
Since $a=\alpha_k^{\frac{2n}{n-1}}$, so $\sqrt{a} = \alpha_k^{\frac{n}{n-1}},$
we see if $k$ is odd this expression is equal to $-2\sqrt{a}$ and if $k$ is even it is equal to $2\sqrt{a}$. Thus, if $k$ is odd, 
$\Rna(b-2\sqrt{a}) = b-2\sqrt{a}$
and if $k$ is even, 
$\Rna(b-2\sqrt{a}) = b+2\sqrt{a}.$
\end{proof}

Thus the $n-1$ points with $k\in \{1,2,3,...,2n-1\},$ for $k$ even we propose are inside of  topological disk hyperbolic components, and the points with $k$ odd are proposed centers of baby $\cM$'s, except in the case of $k=1,2n-1$, where the centers are so near $0$, the components appear to merge into something more complicated. Hence there appear only to be $n-2$ baby $\cM$'s.

\section{Baby Mandelbrots in the $\Rna$ Parameter Space}
\label{sec:babyMandels}

The locations of some baby Mandelbrot sets in the family $R_{n,a,b}$ were established in both \cite{BoydHoeppner1} and \cite{BoydMitchell}. In the present article we are more interested in the newer phenomena of baby Julia sets in parameter space. So, rather than applying similar techniques to establish the location of baby $\cM$'s in this subfamily, we stop at providing some suggestions on how an interested reader could complete this task, by building on the previous articles and the prior results of this article, at least under the assumption that $n$ is sufficiently large. 

By Lemma~\ref{lem:akdefn},
there are $2n-1$ potential centers of hyperbolic components, $n-2$ values of the form $a_k>0,  k=1,2,\ldots,2n-1,$ with $v_-$ fixed for $k$ odd, and $n-1$ values with $v_-$ mapping to $v_+$ for $k$ even.

\begin{conjecture}\label{conj:babyMs}
    For $n\geq 3$, the $n-2$ parameters $a_k$ for  odd $k=3,\ldots, 2n-3$, are each a center of a baby Mandelbrot set in the $a$-plane for the family $r_{n,a}.$ (Note $a_1$ and $a_{2n-1}$ are excluded.)
\end{conjecture}

To prove this conjecture using Douady and Hubbard's Theorem~\ref{thm:DH3}, one must define regions $W_k$ in the $a$-plane, which are closed topological disks with $a_k\in W_k$,
and show that (1): for any $a$ in one of these $W_k$'s, the map $\Rna$ restricted to  $U'_{a,k}$ has a critical point in the $U'_{a,k}$, and the map is polynomial-like of degree 2 on the orbits that remain bounded in $U'_{a,k}$. So, one needs that $\Rna$ is a degree 2 proper, analytic map of each $U'$ to $U^{-}$ where $U'$ is contained in $U^{-}$, with $U'$ relatively compact in $U^-$, because by Proposition~\ref{prop:Uellipse}, each $U'_{a,k}$ maps 2:1 onto to $U^-_a$ for $k$ odd.

Also, one must show (2): as $a$ loops around $\partial W_k$, $v_-$ loops around $U^- \setminus U'_{a,k}$.

Now, recall $a_k = \left(  ({1-e^{i\frac{k\pi}{n}}})/{4} \right)^{\frac{2n}{n-1}}$.
Note we're excluding $a_1$ and $a_{2n-1}$, as those are the parameters that tend to $0$ as $n\to\infty$, since $(2n)/(n-1)\to 2$ and $e^{i\pi (2n-1)/n} = e^{-i\pi /n},$ so $(1-e^{\pm i\pi/ n})\to 0^{-}$ as $n\to\infty$. Excluding them one has a chance to show the other $W_k$'s are bounded away from $a=0$ (and the negative real axis), with a bound depending on $n$. Computer generated images show that there don't tend to be distinct hyperbolic components for $a_1$, and $a_{2n-1}$, but rather that the components that would correspond to potential baby $\cM$s for $a_1$ and $a_{2n-1}$ have ``collided'', with the origin trapped between them. If this is the case, then the two potential baby $\cM$s cannot be baby $\cM$s, and the component where $D_-=D_+$ near $a=0$ seems to be deformed into one unusual component. See Figure~\ref{fig:FixedCritPointExamples1}.

One could initially aim to establish this conjecture at least for all $n$ sufficiently large, or even just start with the $W$s ``farthest'' from the origin (where $a=0$) and negative real axis, which is the branch cut of $a^{1/2n}$; that is, one could start with $W_n$ if $n$ is odd, which should be roughly centered about the positive real axis, or with $W_{n\pm 1},$ if $n$ is even, which are located in symmetrical positions above and below the positive real axis.

One approach to defining the $W_k$'s would be to apply the same idea from \cite{BoydHoeppner1} and define regions $W_k$ that should contain baby $\cM$s in the $a$-plane by using four curves suggested by the definition of $U'_{a,k}$. That would yield the following.

\begin{defn} \label{defn:Wkodd} 
Let $a_k$ be the parameters defined in Proposition~\ref{prop:c-patterns}, so at $a_k$, $v_-=w_k=|a|^{1/2n} e^{i({\psi+2k\pi})/{2n}}$ is a non-canonical critical point, which $r_{n,a_k}$ fixes for $k$ odd,  and for $k$ even maps to $v_+$.
Define $W_k$ for each $k$ as the region in the $a$-plane bounded by the implicit curves:
\begin{align*}
\beta = \left \{  a =
\frac{1}{16}
\left(  \frac{|a|^{1/n}}{2}{e^{i\theta}}- a^{1/2n} \right)^2 
\ \ | \ \ 0 \leq \theta \leq  2\pi \right \}
\\
\tau = \left \{ a= \frac{1}{16} \left( 2{e^{i\theta}} - a^{1/2n} \right)^2 
\ \  | \ \  0 \leq \theta \leq  2\pi \right \}
\\
\rho_k^+= \left \{ a 
= \frac{1}{16}\left( {x}e^{i\left(\Arg(w_k)+\frac{\pi}{2n}\right)}  - a^{1/2n} \right)^2 
\ \  | \ \  0 < x \leq  2 \right \}
\\ 
\rho_k^- = \left \{ a 
= \frac{1}{16}\left( {x}e^{i\left(\Arg(w_k)-\frac{\pi}{2n}\right)} - a^{1/2n} \right)^2 
\ \ | \ \ 0 <  x \leq  2 \right \}
\end{align*}
where $a^{1/2n}$ means the principal root. \end{defn}

We can define all of these $W_k$'s, though only the ones with $k$ odd are conjectural regions containing baby $\cM$'s. 

 Since $v_-=a^{1/2n}-4\sqrt{a}$, we see $\beta$ is the locus of $a$ for which $|v_-| = |a|^{1/n}/2$,
$\tau$ is where $|v_-|=2,$
and $\rho_k^{\pm}$ is where $\Arg(v_-) = \Arg(w_k) \pm \pi/2n.$
Now comparing with Definition~\ref{defn:U_ak_prime}, it's clear that as $a$ loops around $\partial W_k$, $v_-$ loops around $\partial U'_{a,k}$. Hence (2) comes essentially for free from the definition of the regions $W_k$. 

On the other hand, the region definitions given above are implicit, as the equations have an $a$ on each side, which complicates some of the remaining criteria. However, at least for $n$ large, $|a^{1/2n}|$ is close to $1$, and $1$ was in that spot in the definitions of the $W_k$'s in \cite{BoydHoeppner1}, so this may not be too difficult of a complication.

A place where this implicit definition may complicate the process is in the requirement that each $W_k$ be a (closed) topological disk, as one would need to check that was true, at least for $n$ sufficiently large. 

An alternative approach would be to make the $W_k$'s slightly larger than these implicitly defined sets, but one must then check they are not so large that $v_-$ on the boundary escapes $U_{a,k}$. 
That would follow the approach of \cite{BoydMitchell}, which doesn't define the $W$'s so that the critical value exactly traces the boundary of $U'$, but rather defines the $W$'s as polar rectangles, then looks at the image of each boundary piece to show that the critical value loops around $U\setminus U'$ as the parameter loops around $\partial W$. 

Also, note by Lemma~\ref{lem:M_in_annulus}, from the ``outer'' argument, we see that in letting $\delta=1$, we get $u=11/16$ and so for all $n$ sufficiently large and for all $a\in \CC\setminus \DD(1/8, u=11/16),$ we get that $|v_-|> 1+\delta=2$. Thus, inside the curve $\tau$ where $|v_-|\leq 2,$ we know our $W_k'$s are in $\DD(1/8,11/16) \subset \DD(0,13/16)$ at least for $n$ sufficiently large.
So, consider the set $\beta$, where $|v_-| = |a|^{1/n}/2$. 
Now if we take the $n$ sufficiently large that the $W_k$'s are in $\DD(0,13/16) \subset \DD(0,1),$ we have
 $|a|<1$ in $\beta$, so $|a|^{1/n}/2 < 1/2$. 

Here is another place where the implicit definition in the boundary defined by $\beta$ is a bit messy. Once could consider replacing the $a^{1/2n}$ in that definition by a $\lambda_n$ which doesn't depend on $a$. One would just need to calculate $\lambda_n$ to check that then one has as $a$ traces the lower boundary of $W_k$, $v_-$ traces not $U'_{a,k}$ but $V_{a,k} = \{ \lambda_n \leq |z| \leq 2, | \Arg{a} - \Arg{w_k} | \leq \pi/2n \},$ and one would need $U^-_a \supset V_{a,k} \supset U'_{a,k}$ so that $v_-$ is in $U\setminus U'$ for $a$ in $\partial W_k$. So one would need to calculate $\lambda_n$ carefully so that $\lambda_n \leq |a|^{1/n}/2$ but not so small that $V_{a,k}$ sticks out of $U^-_a$. 

Note from, Lemma~\ref{lem:outerboundaryW} given the upper bound  $|v_-|=2$ in these $W_k$'s we'd have $U'_{a,k} \subset \EE$ for all $a$ in $W_k$. 

We provide all of these suggestions in hopes that someone would take on this project and prove the conjecture. If you do, please contact the first author as we would be very interested in seeing your work.

\section{Topological Disks in the $\Rna$ Parameter Space}
\label{sec:paramdisks}

In this section,  we prove Proposition~\ref{prop: preimage classes} describing the dynamics of maps in the topological disk hyperbolic components, then we prove Theorem~\ref{thm:hypcompdisk}.
We close with final suggestions for further study.

Here, we consider $k=2j,$ for $j \in \{ 1,\ldots, n-1 \}$; that is, $k$ even and $k\in \{ 2, \ldots, 2n-2\}$. 

First, for $a=a_{2j}$ 
as defined in Lemma~\ref{lem:akdefn},
 we have $\mapRna(v_+) = v_+$ and $\mapRna(v_-) = v_+$. Further, $v_-$ itself is a critical point, with
$$v_- = w_{2j} = |a|^{1/2n}e^{i(\frac{\psi + 2(2j)\pi}{2n})}
= |a|^{1/2n}e^{i(\frac{\psi}{2n} + \frac{2j\pi}{n})}
,$$ so that it is one of the non-principle roots of $a^{1/2n}$.

$\cH_{2j}$ is the hyperbolic component in the parameter space which contains $a_{2j}.$ Since $a_{2j}$ is the parameter with $\mapRna(v_-)=v_+$, we call $a_{2j}$ the ``center'' of $\cH_{2j}$, and consider $a_{2j}$ and other $a$'s in $\cH_{2j}$.

\subsection{Dynamical plane}
Now, recall from Lemma~\ref{lem:Ra-polynlikeD+} that $\Dplus$  (or $\Dplus(a)$ if this is needed for clarity) is the Fatou component containing $v_+$ (aka $v_+(a)$) and on $\overline{\Dplus},$ the map $\mapRna$ is conjugate to $p_0(z)=z^2$ via $\phi_a$; not just for $a=a_{2j}$, but for all $a \in \cH = \{ a | \ r_{n,a}$ is hyperbolic $\}$. 

Next, for any $a\in \cH,$ $D_-$ is the Fatou component containing $v_-.$ So for $a$ near $a_{2j}$, $\mapRna(D_-) = D_+$. 

We will assume for $a=a_{2j}$ that $v_- \notin D_+$; that is, $D_- \neq D_+$. 

First we discuss the dynamics at $a=a_{2j}$. 
Now, $v_-$ is not a fixed point for $a=a_{2j},$ and so there are $n$ critical points which map to $v_-$ and are all distinct from it. From Lemma~\ref{lem:critptsmapping}, the critical points  $w_{2\ell-1} = |a|^{1/2n}e^{i(\frac{\psi + 2(2\ell-1)\pi}{2n})}$ for $\ell=1,2,3,...,n,$ map to $v_-$, 
and the critical points 
$w_{2\ell}=|a|^{1/2n}e^{i(\frac{\psi + 2(2\ell)\pi}{2n})} $ for
$\ell=0,1,2,...,n-1,$ map to $v_+,$ with $w_0=v_+$ in this case and our special $j$ satisfying $w_{2j}=v_-$. So,  not counting $v_+$ itself, there are $n-1$ critical points mapping to $v_+,$ one of which is $v_-$, and $n$ critical points mapping to $v_-$.  

For each $\ell \in \{ 0, \ldots, 2n-1 \},$ let $D_\ell$ be the Fatou component containing the critical point $w_\ell$, so $D_0=D_+$ and for $a=a_{2j}$, $D_{2j}=D_-.$

Thus $\mapRna$ maps each $D_{2\ell-1}$
onto $D_-$ as a $2:1$ branched covering ramified over the critical point. 
Similarly, $\mapRna$ maps each $D_{2\ell},$
including $D_-$, 
as a $2:1$ branched cover onto $D_+$,  ramified over the critical point in $D_{2\ell}$.

\begin{figure}
\centering
\begin{subfigure}{.5\textwidth}
  \centering
\includegraphics[width=.95\textwidth,keepaspectratio]{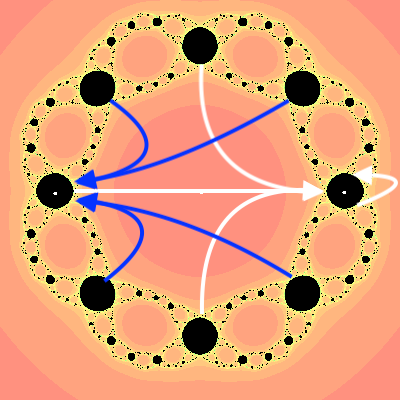}
  \caption{Julia set $n=4, a\approx 0.16,$ }
  \label{Jfixedcrit_n4mapping}
\end{subfigure}
\begin{subfigure}{.5\textwidth}
  \centering
\includegraphics[width=.95\textwidth,keepaspectratio]{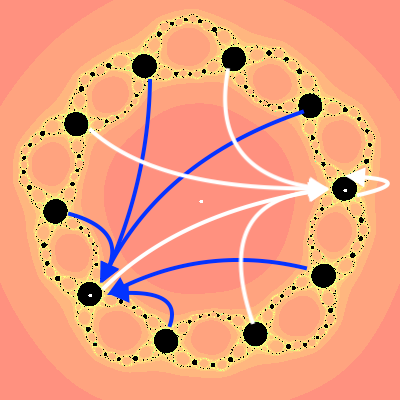}
  \caption{Julia set, $n=5, a\approx 0.11+0.11i$}
  \label{Jfixedcrit_n5mapping}
\end{subfigure}
\caption{Arrows illustrate $\mapRna$'s action on the critical components of the filled Julia set, in the case $\mapRna(v_-)=v_+,$ so $v_-$ is near the center of $D_-$.}
\label{fig:FixedCritPointJuliaMappings}
\end{figure}

See Figure~\ref{fig:FixedCritPointJuliaMappings} for an illustration of how the Fatou components containing critical points map, in these cases of $a=a_{2j}$ so that $v_+$ is fixed with $v_-$ mapping onto $v_+$.

Next, if we consider $a$ near $a_{2j}$, of course $D_+$ varies with $a$, but as shown in Lemma~\ref{lem:Ra-polynlikeD+}, $v_+$ is safely inside of $D_+$. Hence for $a$ near $a_{2j}, v_-$ is safely inside of $D_-$, since $\mapRna(v_-)=v_+$ at $a=a_{2j}$. As $a$ varies away from $a_{2j}, v_-$ moves toward the boundary of $D_-$, its image $\mapRna(v_-)$ moves away from $v_+$, heading toward the boundary of $D_+$ as $a$ grows. Additionally, the critical point preimages of $v_-$, the $w_{2\ell-1}$'s, approach the boundaries of their Fatou components $D_{2\ell-1}$, or, more accurately, the boundary of the $D_{2\ell-1}$'s tend to the $w_{2\ell-1}$'s. Since $\mapRna$ is a $2:1$ ramified cover from each component $D_{2\ell-1}$ of $\mapRna^{-1}(D_-)$ onto $D_-$, as $v_- \to \partial D_-$, when $v_-$ hits the boundary of $D_-$, each $D_{2\ell-1}$ devolves into a topological lemniscate, with a figure 8 boundary with the critical point $w_{2\ell-1}$ at the crossing. 

As $a$ continues to move away from $a_{2j},$ just past the lemniscate phase, there are then two Fatou components that replace $D_{2\ell-1},$ neither containing a critical point, each mapping $1:1$ onto the preimage of $D_+$ that is no longer $D_-$ because $v_-$ has moved out of it.

\begin{figure}
  \centering
  \begin{subfigure}{.45\textwidth}
   \centering
\includegraphics[width=1.0\textwidth,keepaspectratio]{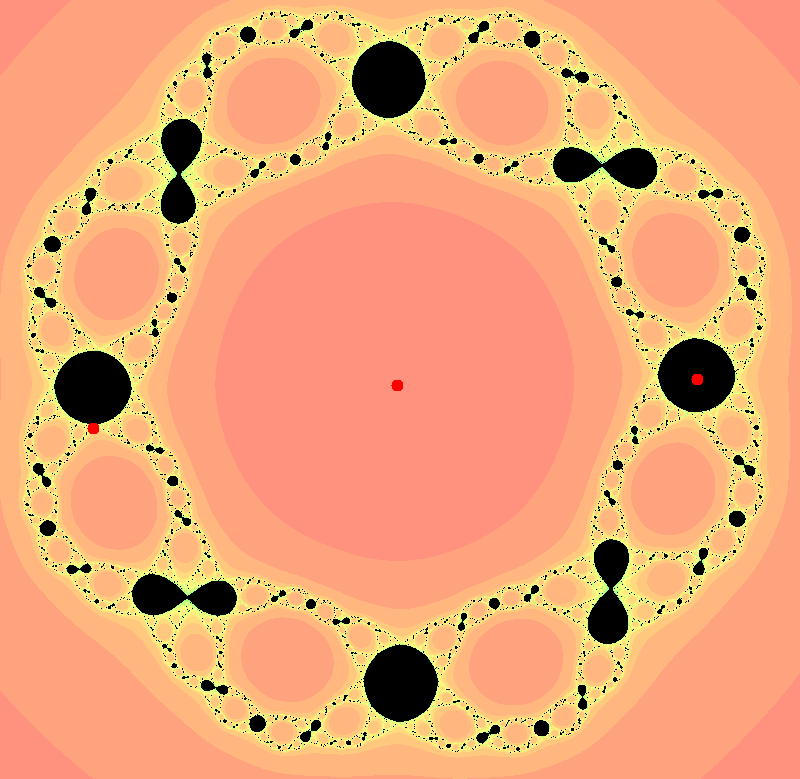}
  \caption{$a=0.16+0.026i$ has $v_-\in\partial D_- $, lemniscates in $J$ }
  \label{Jfixedcrit_lemniscate}
\end{subfigure}
\begin{subfigure}{.45\textwidth}
  \centering
\includegraphics[width=1.0\textwidth,keepaspectratio]{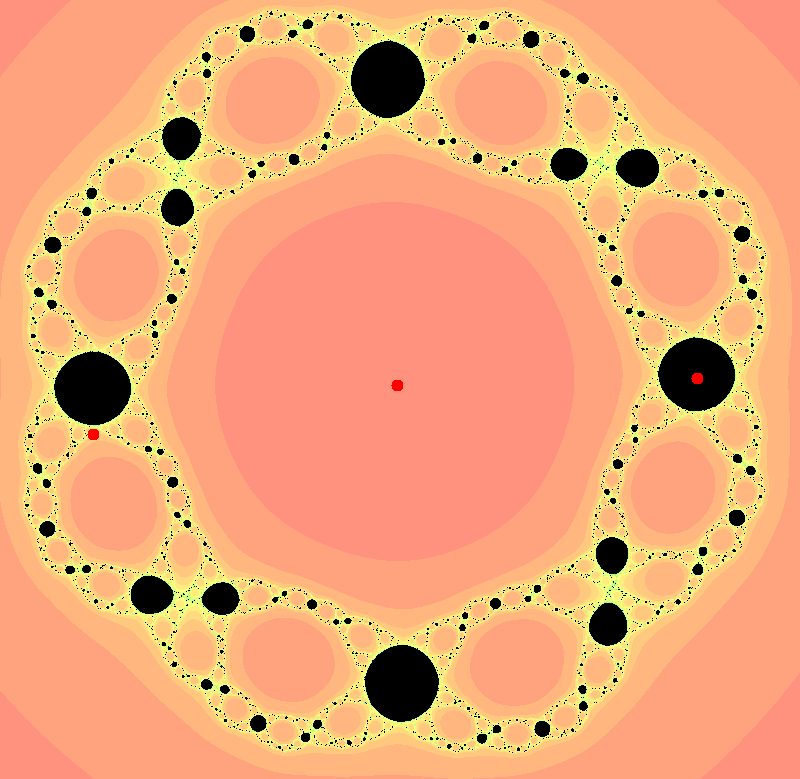}
  \caption{$a=0.16+0.03i$ has $v_-$ just outside of $D_-$}
  \label{Jfixedcrit_twodisks}
\end{subfigure}
\caption{Julia sets for $\Rna, n=4$ showing bifurcation as $v_-$ leaves $D_-$. $v_{\pm}, 0$ are marked in red, with $v_-$ on the left. 
}
\label{fig:lemniscatebifurcation}
\end{figure}

See Figure~\ref{fig:lemniscatebifurcation} for a Julia set with $v_- \in \partial D_-$ and another with $v_-$ just outside of $\overline{D_-}.$

To describe further the combinatorial dynamics of these types of maps,
we first apply Theorem~\ref{thm: angle assignments} to assign sets of external angles to the subset of the Julia set consisting of $\partial \Dplus$ and its tree of preimages. Note that this result is written to apply to any polynomial-like generalized McMullen map $R_{n,a,b}$ whose restriction to the filled Julia set within $U'$ is conjugate to the restriction of a quadratic polynomial $P_c$ on its filled Julia set.

We apply this result here where $\overline{D_+}$ is the baby filled Julia set $K_+$, and the quadratic polynomial we are conjugate to is $P_0(z)=z^2$ with filled Julia set $\overline{\DD}$. Then we get a set of angle assignments on each preimage of $\partial D_+$, which respects the dynamics of $r_{n,a}$.

In \cite{BoydBrouwer1}, we applied that theorem to describe combinatorially various preimages of $K_+$ in terms of external angle identifications; e.g., in ``the basilica'' map $P_{-1}(z) = z^2-1$, angles $1/3$ and $2/3$ meet at a fixed point of $P_{-1}$. In the present case, for $\Rna$, the baby Julia set $\overline{D_+}$ is always a quasi-disk, so no angles are identified. But we can still give some description in terms of these angle assignments that sheds light on the dynamics.

\begin{prop} \label{prop: preimage classes}

    Let $r=r_{n,a}$ be in the subfamily given in~(\ref{eqn:subfamilyrna}),
    so the Julia set $J(r)=J_r$ contains a baby Julia set $\overline{D_+}$ on which $r$ is conjugate to $P_0(z)=z^2$ on its Julia set, the closed unit disk, via a map $\phi_a$ which is q.c.\ in a neighborhood of the filled Julia sets and analytic on the interior.
     Assume in addition that $v_-$ lies in a preimage of $D_+$.

    Then $J_r$ also contains the infinite tree of preimages $J_* = \displaystyle \cup_{m=0}^\infty\  r^{-m}(\partial D_+)$, where each preimage component satisfies one of the following:
    \begin{enumerate}
        \item For components $D_{2\ell}$, $\ell=0,\dots,n-1$, each containing the critical point $w_{2\ell}$, we have  
       \\
        $r\Big(\Gamma|^{\partial D_{2\ell}}\big([0,\frac12)\big)\Big) = r\Big(\Gamma|^{\partial D_{2\ell}}\big([\frac12,1)\big)\Big) = \Gamma|^{\partial D_+}\big([0,1)\big)$.
        This includes $\Dplus = D_{0}$, and $\Dminus$, for some $\ell>0$.

        \item  For components $\Dmke$ of $r^{-m}(D_{2\ell})$ for $\ell=1,\dots,n-1$ and $m\geq1$, we have \\
        $r^{m+1}\Big(\Gamma|^{\partial \Dmke}\big([0,\frac12)\big)\Big) = r\Big(\Gamma|^{\partial \Dmke}\big([\frac12,1)\big)\Big) = \Gamma|^{\partial D_+}\big([0,1)\big)$.
        
        \item For components $D_{2\ell-1}$, $\ell=1,\ldots,n,$ each containing the  critical point $w_{2\ell-1}$, we have \\
        $r^2\Big(\Gamma|^{\partial D_{2\ell-1}} \big([x,x+\frac14)\big)\Big) = \Gamma|^{\partial \Dplus}\big([0,1)\big)$, for $x=0,\frac14,\frac12,\frac34$.

        \item For components $\Dmko$ of $r^{-m}(D_{2\ell-1})$ for $\ell=1,\dots,n$ and $m\geq1$, we have \\
        $r^{m+1}\Big(\Gamma|^{\partial \Dmko} \big([x,x+\frac14)\big)\Big) = \Gamma|^{\partial \Dplus}\big([0,1)\big)$, for $x=0,\frac14,\frac12,\frac34$.
    \end{enumerate}
\end{prop}

\begin{proof}
First recall the statements about how $r$ maps each component as described just above Theorem~\ref{thm: angle assignments}. 

For (1), note $r(D_{2\ell})=\Dplus$ for each $\ell=0,\dots,n-1$, where each  $r(D_{2\ell})=\Dplus$ is a $2:1$ mapping. 

For (3), note for each $\ell=1,
        \dots,n$, $r(D_{2\ell-1})=\Dminus$ as a $2:1$ mapping, so $r^2(D_{2\ell-1})=\Dplus$ is a $4:1$ mapping.

For (2), note as all critical point components have already been identified, we must have that $\Dmke$ maps $1:1$ onto its image, and $r^m(\Dmke)=D_{2\ell}$, so $r^{m+1}(\Dmke)=\Dplus$ in a $2:1$ mapping. Note the case of $\ell=0$ is excluded.

For (4), note as noted above, $\Dmko$ cannot contain a critical point and so must map $1:1$ onto its image, where $r^m(\Dmko)=D_{2\ell-1}$. Thus $r^{m+1}(\Dmko)=\Dplus$ in a $4:1$ mapping. 

Then each of the statements about how subintervals $[0,\frac12)$, $[\frac12,1)$, $[0,\frac14)$, $[\frac14,\frac12)$, $[\frac12,\frac34)$, and $[\frac34,1)$ map follows directly from applying Theorem~\ref{thm: angle assignments}. 
\end{proof}

Just before $v_-$ is about to escape $K(r)$, we observe that the preimages of $D_-$ had boundaries that are nearly figure 8's. In \cite{BoydBrouwer1}, we detected ``altered'' baby Julia sets which were preimage copies of a true baby Julia set with some angle identifications different from what one would see in a quadratic polynomial. Again, in our case of interest there are no angle identification changes, but we do have $J(r)$ which consists of a tree of preimages of the quasi-disk $D_+$ that look one way and the preimages of $D_-$ which have the two lobes of the (near) figure 8 each mapping onto $D_-$, and then their preimages are ``clone'' conformal copies of figure 8's mapping onto figure 8's.  So there are alterations in the shape in a sense, though they're still all quasi-disks. 

Note that in this subfamily, when $a$ is varied just enough for $v_- \notin K(r)$, so the assumption that $v_-$ maps into $D_+$ is a near miss---like in the right side of Figure~\ref{fig:lemniscatebifurcation}---some of these statements still apply. We still have item (1), except that $D_-$ does not exist per se--there is still a preimage of $D_+$ near $v_-$, but it does not contain $v_-$. We still have item (2). Items (3) and (4) no longer apply as stated because the odd critical points escape, \textit{but} consider the component that for nearby $a$ was $D_-$; it is still a preimage of $D_+$, so it looks the same, but since $v_-$ is not inside of it, its preimage components, rather than each containing a critical point and mapping $2:1$ onto $D_-$, consist of a pair of components each mapping $1:1$.


\subsection{Hyperbolic components}

The goal of this subsection is to establish Theorem~\ref{thm:hypcompdisk} and provide suggestions on improving the results in the future.

\begin{proof}[Proof of Theorem~\ref{thm:hypcompdisk}]
At the beginning of Section~\ref{sec:paramdisks} we defined $\cH_{2j}$ as the hyperbolic component in the bifurcation locus of $\Rna$ with center $a_{2j}$   such that  
$r(v_+) = v_+$ is a super-attracting fixed point, $r(v_-)=v_+,$ and $v_-$ is a critical point and in particular a non-principle root of $a^{1/2n}$. This definition was given in Lemma~\ref{lem:akdefn} as to match the statement of Theorem~\ref{thm:hypcompdisk}. 

We know the behavior at each $a_{2j}$ is different, so each $a_{2j}$ has to be in a different hyperbolic component, excluding the degenerate cases which are in the same component near $a=0$, as described at the beginning of Section~\ref{sec:dpr}: those components would correspond to $j=0,n+1$ which is why the statement of this theorem listed only $j=1,\ldots,n$.

What remains is to define and describe the map $\Phi_j$, for $j=1,\ldots,n$.

Recall from Lemma~\ref{lem:Ra-polynlikeD+} that for $a=a_{2j},$ which is safely interior to $\cH_{2j},$ $\mapRna$ is conjugate to $p_0(z)=z^2$ via a map $\phi_a$ which quasi-conformally maps a neighborhood of $\overline{D^+}$ to a neighborhood of $\overline{\mathbb{D}}$, and $\phi_a$ from the filled Julia set to the closed disk is conformal.

Note by Theorem~\ref{thm:DHC2T1}, $\phi_a$ given by Lemma~\ref{lem:Ra-polynlikeD+} depends continuously on $a$ for $a$ in the interior of any hyperbolic component since the baby Julia set $\Dplus$ is always connected.

Define the map $\Phi_j \colon \cH_{2j} \to \DD$ by
\begin{equation}
\label{eqn:topdiskHmap}
    \Phi_j(a) = \phi_a (\Rna(v_-(a))), 
\end{equation}
remembering $v_-$ is a function of $a$. 
That is, $\Phi_j(a) = \phi_a ( \Rna ( \Rna ( w  )  )  )$, where $w \in \mapRna^{-1}(v_-)$, i.e.,  $w^{2n}=a$ and $w$ is any of the $n$ critical points that map to $v_-$.  Specifically, 
from Lemma~\ref{lem:critptsmapping}, at $a=a_{2j}$ the center of $\cH_{2j}$, the $n$ critical points  $w_{2l-1} = |a|^{1/2n}e^{i(\frac{\psi + 2(2l-1)\pi}{2n})}$ for $l=1,2,3,...,n,$ map to $v_-$. 
For the purpose of demonstration, we use $l=1$, remembering that no value of $l$ makes this point equal to $v_-$. 
Set $f_{n}(a) = |a|^{1/2n}e^{i(\frac{\psi + 2\pi}{2n})}=|a|^{1/2n}e^{i(\frac{\psi}{2n} + \frac{\pi}{n})}=w_1$, so $f_{n}(a)$ is the map sending $a$ to a fixed non-principle root of $a^{1/2n},$ the one that is argument $\pi/n$ larger than $v_+$.

So in the hyperbolic component $\cH_{2j},$ 
we have 
$\Phi_j =\phi_a \circ \Rna \circ \Rna \circ f_{n},$
where $f_n(a)$ is a consistent choice of $a^{1/2n},$
$\mapRna$ is a rational map in both applications a $2:1$ branched covering ramified over a critical point, and we know $\phi_a$ is analytic in $z$ on $\overline{D_+(a)}$ by Lemma~\ref{lem:Ra-polynlikeD+}. 

Now to shed light on the map $\Phi_j$, we'll view it a different way. 
Since $a_{2j}$ is the dynamical center of a hyperbolic component which is where $v_-$ maps to the fixed critical value $v_+$, let $\hat{a}_{2j}$ be a point near $a_{2j}$ in the same component, so that its map has no non-persistent critical orbit relation since now we have that $r_{n,\hat{a}_{2j}}(v_-(\hat{a}_{2j}))\neq v_+(\hat{a}_{2j})$. Thus,
$r_{n, \hat{a}_{2j}}$
is topologically stable, so by McMullen and Sullivan \cite{McMSull1998} (as described in Section~\ref{sec:prelim}), we can use this as a basepoint to define a holomorphic motion of the entire sphere, and as there are no obstructions, this holomorphic motion extends to the punctured hyperbolic component, having removed its dynamical center. This holomorphic motion respects, or, more specifically, conjugates, the dynamics; i.e,
 for $\cH_{2j}^*:= \cH_{2j}\setminus \{ a_{2j}\}$, 
 let $\tau^{2j} : \cH_{2j}^* \times \Chat \to \Chat$ be this holomorphic motion, with notation $\tau^{2j}(a,z) = \tau^{{2j}}_a (z)$, 
 and such that $\tau^{2j}({\hat{a}_{2j}},\cdot) = id,$ and 
 \begin{equation}\label{eqn:tau-r}
 \tau^{2j}_a \circ r_{n, \hat{a}_{2j}} = r_{n,a} \circ \tau^{2j}_a.     
 \end{equation}
Now, since $\tau$ is a holomorphic motion, for each fixed $a\in \cH_{2j}^*$, $\tau^{2j}_a$ is an injective function of $z$ which again is the identity at the basepoint, and for each fixed $z$ in the dynamical plane of the basepoint map, $\tau^{2j}(\cdot, z)$ is a holomorphic function of $a$.

Due to uniqueness from conjugacy with the dynamics, we can rewrite the definition of $\Phi_j$ using this $\tau$. Our original map is $\tau^{2j}: \cH_{2j}^* \times \Chat \to \Chat$, but consider the restriction to  $\cH_{2j}^* \times \overline{\Dplus}(\hat{a}_{2j}) \to \overline{\Dplus}(a)$; that is, let us restrict to the baby filled Julia set at the reference point $\hat{a}_{2j}$ mapping to the baby filled Julia set at the target $a$. If you fix an $a\in \cH_{2j}^*$ then this restricted $\tau$ maps $\overline{\Dplus}(\hat{a}_{2j})$ to $\overline{\Dplus}(a)$. See the left half of Figure~\ref{fig:tau-r-phi}. 
Now, the map $\phi_a$ maps $\Dplus(a)$ to $\DD$, so if you then apply $\phi_a$ you again get the topological disk $\overline{\DD}$. 
Instead of as above in the definition of $\Phi$ just considering where each different $r(v_-)$ is in its own $D_+(a)$ before it is mapped into $\DD$ by $\phi_a$, which depends on where $v_-(a)$ is in $D_-(a)$, 
think about where $r_{\hat{a}_{2j}}(v_-(\hat{a}_{2j}))$ or its preimage $v_-(\hat{a}_{2j})$, is moved to via $\tau$ into each $D_+(a)$ or $D_-(a)$, respectively. This ends up being equivalent thinking, due to the uniqueness of conjugation with dynamics guaranteed by McMullen-Sullivan.

More precisely, by Equation~\ref{eqn:tau-r}, and since $\tau_{\hat{a}_{2j}}^{2j}=id$, we have
$$
\tau_{a}^{2j}(r_{n,\hat{a}_{2j}}(v_-(\hat{a}_{2j})) 
= r_{n,a}(v_-(a)) \in D_+(a).
$$
Now if we apply $\phi_a$, we map to a point in $\DD$, as it is the Julia set of $z\mapsto z^2$.
Since $a\neq a_{2j}$, $r_{n,a}(v_-(a)) \neq v_+(a)$, so we know $\phi_a(r_{n,a}(v_-(a))) 
\neq 0$. 
So, $$
\Phi_j(a) = \phi_a (r_{n,a}(v_-(a))) = \phi_a(\tau_{a}^{2j}(r_{n,\hat{a}_{2j}}(v_-(\hat{a}_{2j})) ).
$$
See Figure~\ref{fig:tau-r-phi}.

\begin{figure}  
  \centering
\includegraphics[width=1.0\textwidth,keepaspectratio]{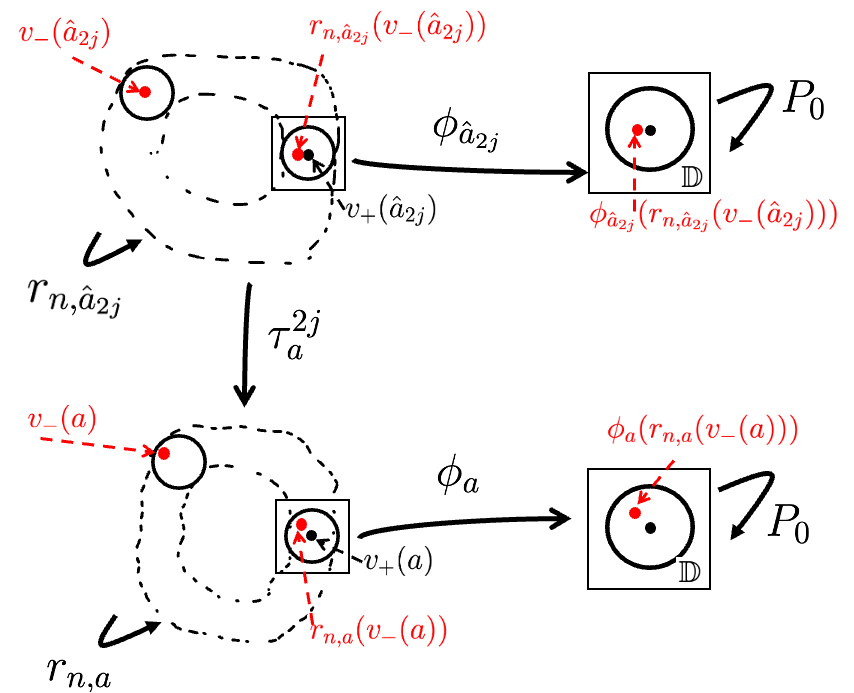}
\caption{\label{fig:tau-r-phi}
Diagram illustrating the maps $\tau, \phi, \Phi, r$ and $p_0$, and the Julia sets for $\Rna$  (left) and $P_0$ (right). The point of interest $v_-(a)$ and its images under various maps are marked in red; $v_+$ (left) and $0$ (right) are shown as black dots.
}
\end{figure}

Remember in viewing the diagram of Figure~\ref{fig:tau-r-phi} that $r_{n,\hat{a}_{2j}}$ is just one map, and in the above picture $r_{n,\hat{a}_{2j}}(v_-)$ stays put in its $D_+$, then as $a$ varies below, its image under $\tau$, which is $r_{n,a}(v_-(a))$, varies. Again thinking of the top half of the picture, $\Dplus(\hat{a}_{2j})$ is a quasi-disk as shown in Lemma~\ref{lem:Ra-polynlikeD+}, and for $r_{n,\hat{a}_{2j}},$ one critical point, $v_+$, is fixed, and the other, $v_-$, must lie in a pre-image of $D_+$, because $\hat{a}_{2j}$ is near $a_{2j}$ in its hyperbolic component, and at $a_{2j}$, $v_-$ maps to $v_+$.  So $D_-(\hat{a}_{2j})$ is also a quasi-disk.  

Now, as $a$ varies in $\cH_{2j}^*$, $\tau_a^{2j}$ is a holomorphic function of $a$ (and injective in $z$), and $\phi_a$ depends continuously on $a$ for any $a$ in the interior of a hyperbolic component by Lemma~\ref{lem:Ra-polynlikeD+}. We argue that $\Phi(a)$ sweeps out $\DD\setminus \{ 0 \}$ as $a$ sweeps out $\cH_{2j}^*$.
We have that $v_-$ is essentially $a^{1/2n}-4\sqrt{a}$. 
One theoretical problem for $\Phi_j$ could be 
a problem with a square root or $1/2n$ root, which can happen for hyperbolic components which run into $a=0$, but we have excluded these with our restrictions on $j$. Or, another problem would be a critical orbit interfering, but we are controlling this with our hypotheses on this subfamily that one critical orbit is fixed and the other is the only one which is independent. 

Essentially, since there are only two critical values, $v_+$ is stuck in a hyperbolic component avoiding $a=0$, so $v_-$ can have its ``full experience'' as the parameter $a$ varies in $\CC$ in a neighborhood of $a_{2j}$. We have a quasi-disk $\Dplus$ in dynamical space, and which has a preimage $\Dminus$ on which $\mapRna$ is a 2:1 branched covering to $\Dplus$, so $D_-$ is a quasi-disk. Since $r(v_-)$ can be ``anywhere'' in $D_+$ as one varies $a$, we see that $v_-$ can be anywhere in $D_-$, so, $\phi_a(\mapRna(v_-))$ sweeps out the whole unit disk as $a$ varies in $\cH_{2j}$. Now we define ``anywhere'' as relative to the corresponding points for the map $r_{n,\hat{a}_{2j}}$.  That is, as $a$ varies in a neighborhood about $\hat{a}_{2j}$, a point sweeps out the $\DD$ (in the upper right of the figure) which is the image under $\phi_{\hat{a}_{2j}}$ of $D_+(\hat{a}_{2j})$, thus the preimage under $\phi_{\hat{a}_{2j}}$ sweeps out the disk $D_+(\hat{a}_{2j})$, hence the preimage under $r_{n,\hat{a}_{2j}}$ sweeps out the disk $D_-(\hat{a}_{2j})$, and that is the reason you see a disk in parameter space.
That gives $\Phi$ is surjective, and its injective since $\tau$ is a holomorphic motion and $\phi$ is a q.c.-homeomorphism, and it's continuous because $\tau$ and $\phi$ vary continuously with $a$. This yields a copy of $\DD$ in parameter space, since $v_-$ is just calculated from $a$ by simple complex arithmetic and roots, and our hyperbolic components avoid $0$.

Thus $\Phi_j$ is a homeomorphic map of the hyperbolic component $\cH_{2j}$ to the unit disk $\DD$. This is why we see topological disks in the parameter space.
\end{proof}

Finally, now that we have established Theorem~\ref{thm:hypcompdisk}, in the next subsection, we point the way for its improvement.

\subsection{Future work on the location of the topological disk components}
To establish more specific results about these topological disk components for $\Rna$, such as defining explicit regions for each $n$ in which they are contained,
one could first attempt to prove Conjecture~\ref{conj:BabyJcriteria} in general then apply it to this family, or vice versa, to establish this conjecture in the simpler case first to point the way to the future work of generalization. 

For this family of maps $\Rna$, one could try to define some regions $\sA_n$ in parameter space, for which one could calculate or bound their distance from $0$ in terms of $n$, since $0$ is to be avoided as it is where $D_-=D_+$. The regions could be defined dynamically or more simply using polar rectangles. The approach would be to show that for the boundary of these regions, or at least for one of these regions as a start, $r(v_-)$ loops around the boundary of $D_+$.

Considering computer generated images such as those provided in Figure~\ref{fig:FixedCritPointExamples1}, one can observe that if $n$ is even there’s a topological disk centered on the positive real axis in dynamical plane, 
in a region we denote $\sA_n$, 
and if $n$ is odd, there are two symmetrical ones in quadrants 1 and 4, which would have the indices ${n \pm 1}$. 
Given this line of thinking,
the ``worst'' case is $n=3$ in which $\sA_2$ and $\sA_4$ are the only disks.
A starting approach could be to start with larger $n$ to make the argument easier, then consider decreasing $n$. 
Perhaps one can show that $\sA_n$ or $\sA_{n \pm 1}$ is outside of a lower bound on radius and then say in the right half plane, although this would not be true for $n=3$, where the cardioid has cusp at $0$ and the two disks look centered along the imaginary axis.  
As of now, we don't know how big the topological disks or baby Mandelbrot's might be, but we do have the centers. 

One could try using Lemma~\ref{lem:M_in_annulus} or Corollary~\ref{cor:M_in_disk}, showing that for any $\ep>0$, there is an $n$ sufficiently large that, if $l \leq 1/32$ and $u> 3/8$, then $M_n(\Rna)\subset \mathbb{A}(l, u)+1/8 \subset \DD(0, 1/2+\ep). $  In proof of the Lemma~\ref{lem:M_in_annulus}, we also showed if $a \in \DD(1/8,1/32)$ (i.e., $|a-1/8| < 1/32$) then $|v_-| \leq L(\infty) < 0.614,$ and if $|a| > 1/2$, then for sufficiently large $n$, $|v_-| \geq 1+\ep_1$, for $\ep_1$ depending on how far $u$ is from $3/8$. 
One could start with an annulus and  try to establish one or two (depending on the parity of $n$) topological disk components to the right of the line $\Re(z)=1/8$. 
The reason for staying to the right of the line is that 
one must avoid $a$ close to $0$, for which $v_-$ too close to $v_+$. One needs the Fatou components of $v_-$ and $v_+$ staying disjoint, to avoid the degenerate case discussed near the beginning of Section~\ref{sec:dpr}.

Furthermore, remember $\Arg$ has range $(-\pi, \pi]$, and $a_k$ is the parameter such that $v_- = w_k$ the $k^{th}$ critical point. Note $w_{2j}$ is not the critical point of $\mapRna$ for any $a \in W_{2j}$, but the specific critical point of $\mapRna$ for the hyperbolic center $a_{2j}$.

A strategy could be to   
 find bounds in terms of $n$ for the location of
 $a$, $v_-$ and/or  $w_{2j}$, for any $a \in \sA_{2j}$ (with center $a_{2j}$ at which $v_- = \crpt_{2j}$), for $j\in \{ 1, \ldots n-1\},$ (so excluding $a_0=0$). 
As a start, from definition of $W_{2j}$,
we know for every $a$ in the interior of $W_{2j}$, we have 
$\frac{|a|^{1/n}}{2} < |v_-| < 2$ 
and 
$ | \Arg(v_-) - \Arg(w_{2j}) | < \frac{\pi}{2n}.  $
Depending on the definition of the $W_k$'s one may need to prove that the $W$'s partition up a region containing the boundedness locus in a useful way.


\bibliographystyle{plain}

\end{document}